\documentstyle[bezier,grcalc,amstex,12pt]{amsart}
\evensidemargin.cm

\def\myauthor{Ph\`ung H{{\accent"5E o}\kern-.28em\raise.2ex\hbox{\char'22}\kern-.20em} H{a\kern-.370em\raise.16ex\hbox{\char'47}\kern.1em}i}

\def\myperaddress{Hanoi Institute of Mathematics, P.O.Box 631, 10000 Boho, Hanoi}
\def\mythanks{The work was done during the author's stay at the Abdus Salam International Center for Theoretical Physics, Trieste, Italy. He would like to thank the Center for its hospitality and financial support.}
\def\myabstract{Given a Hecke symmetry $R$, one can define a matrix bialgebra $E_R$ and a matrix Hopf algebra $H_R$, which are called function rings on the matrix quantum semi-group and matrix quantum groups associated to $R$. We show that for an even Hecke symmetry, the rational representations of the corresponding quantum group are absolutely reducible and that the fusion coefficients of simple representations depend  only on the rank of the Hecke symmetry. Further we compute the quantum rank of simple representations. We also show that the quantum semi-group is ``Zariski'' dense in the quantum group. Finally we give a formula for the integral.}

\def\thedate{September 21, 1998}

\def\amshead{\date{\thedate}
\title{On Matrix Quantum Groups of Type $A_n$ }
\author{Phung Ho Hai}
\address{\myperaddress}
\curraddr{Max-Planck Institut f\"ur Mathematik, Gottfried-Claren-Str. 26, 53225, Bonn, Germany}
\thanks{\mythanks}
\keywords{Matrix quantum groups, Hecke symmetry, Hecke algebras, Murphy operator.}
\email{phung@@mpim-bonn.mpg.de}
\begin{abstract}\myabstract\end{abstract}
\subjclass{Primary 16W30,17B37 , Secondary 17A45, 17A70}
\maketitle }

\let\testint\int
\def\int{{\displaystyle\testint} }
\newcommand{\bK}{{\Bbb K}}
\newcommand{\bZ}{{\Bbb Z}}

\newcommand{\bN}{{\Bbb N}}
\newcommand{\BBb}[1]{{\Bbb #1}}

\newcommand{\Ln}[1]{\overline{L_{#1}}}
\newcommand{\ove}[1]{\overline #1}
\newcommand{\tttt}{\end{equation}}
\def\rank{\mbox{rank}}

\newcommand{\mono}[3]{#1_{#2_1}^{#3_1}#1_{#2_2}^{#3_2}\cdots #1_{#2_n}^{#3_n}}

\newcommand{\edim}{\mbox{\rm 8dim}}
\newcommand{\etr}{\mbox{\rm 8tr}}
\newcommand{\rdim}{\mbox{\rm rdim}}

\def\lora{\longrightarrow}
\def\ot{\otimes}

\def\loma{\longmapsto}
\def\Vn{{V^{\ot n}}}

\newcommand{\bbas}{\begin{eqnarray*}}
\newcommand{\eeas}{\end{eqnarray*}}
\newcommand{\bbar}{\begin{array}}
\newcommand{\eear}{\end{array}}
\newcommand{\bbs}{\begin{displaymath}}
\newcommand{\ees}{\end{displaymath}}
\newcommand{\bb}{\begin{equation}}
\def\ee{\end{equation}}
\def\eea{\end{eqnarray}}

\def\bba{\begin{eqnarray}}

\newtheorem{thm}{Theorem}[subsection]
\newtheorem{lem}[thm]{Lemma}
\newtheorem{dn}[thm]{Definition}
\newtheorem{edl}[thm]{Theorem}
\newtheorem{cor}[thm]{Corollary}
\newtheorem{pro}[thm]{Proposition}

\def\Lambd{{\mit\Lambda}}

\def\det{\mbox{\rm det}}

\def\Im{\mbox{\rm Im}}

\def\Hom{\mbox{\rm Hom}}

\def\End{\mbox{\rm End}}

\def\H{{\cal H}}
\def\W{{\cal W}}

\def\V{{\cal V}}
\def\P{{\cal P}}

\def\eee{ {\mbox{\normalsize\rule[-.25ex]{.75ex}{1.5ex}}\\}}
\def\eqeee{\mbox{\normalsize\rule[-.25ex]{.75ex}{1.5ex}}}

\newcommand{\proof}{{\it Proof.\ }}
\newcommand{\va}{\varepsilon}
\def\nml{\normalsize}

\newcommand{\tr}{\mbox{\rm tr}}
\renewcommand{\dim}{\mbox{\rm dim}}

\def\rref#1{(\ref{#1})}

\font\Fraktur=eufm10 scaled\magstep1          
   \newcommand{\fraktur}[1]{\mbox{\Fraktur #1}}  %
   \font\Fraktu=eufm7 scaled\magstep1            
   \newcommand{\fraktu}[1]{\mbox{\Fraktu #1}}    %
   \font\Frakt=eufm5 scaled\magstep1             
  \newcommand{\frakt}[1]{\mbox{\Frakt #1}}      %
   \def\fr#1{\mathchoice{\fraktur {#1}}            
                        {\fraktur {#1}}            
                        {\fraktu {#1}}             
                        {\frakt {#1}}  }           

\newcommand{\Ss}{\fr S}

\def\db{{\mathchoice{\mbox{\rm db}}
                    {\mbox{\rm db}}
                    {\mbox{\scriptsize\rm db}}
                    {\mbox{\tiny\rm db}} }}

\def\ev{{\mathchoice{\mbox{\rm ev}}
                    {\mbox{\rm ev}}
                    {\mbox{\scriptsize\rm ev}}
                    {\mbox{\tiny\rm ev}} }}

\def\id{{\mathchoice{\mbox{\rm id}}
                    {\mbox{\rm id}}
                    {\mbox{\scriptsize\rm id}}
                    {\mbox{\tiny\rm id}} }}

\newcommand{\Comod}{\mbox{Comod}}

\def\part{\vdash}

\def\lam{{\lambda}}
\def\unit{{\mathbf 1}}

\begin{document}

\bibliographystyle{plain}
\amshead

\section{Introduction}\label{sectinto}

Matrix quantum groups of type $A_n$ generalize the standard deformations of the matrix group GL$(n)$ introduced by Faddeev, Reshetikhin, Takhtajian \cite{frt}. They are defined in terms of Hecke symmetries. To a Hecke symmetry $R$ there are associated a bialgebra $E_R$ and a Hopf algebra $H_R$, which are considered as the ring of regular functions on the corresponding matrix quantum  semi-group and group of type $A$, respectively. If $R$ is the Drinfel'd-Jimbo solution to the Yang-Baxter equation of type $A_n$,  the Hopf algebra $H_R$ reduces to the above mentioned deformations. 

Quantum groups of type $A$ were firstly studied by Gurevich \cite{gur}, who generalized\hyphenation{Lyu-ba-shen-ko} Lyubashenko's results on vector symmetry for the case of Hecke symmetry. Many interesting properties of quantum matrix (semi-) groups were found \cite{gur1,ls,ph97,ps96,gur96}. These results show on one hand the similarity in the theory of quantum matrix groups of type $A_n$ and the classical groups  GL$(n)$, and on the other hand the relations between quantum matrix groups of type $A_n$ and the Hecke algebras.

In this work we study rational representations of matrix quantum groups of type $A_n$, i.e., comodules of Hopf algebras associated to even Hecke symmetries. For an even Hecke symmetry, one has a notion of rank, which is, in the classical case, equal to the dimension of the vector space, on which the Hecke symmetry is defined. It was showd by Gurevich that, in general, the rank and the dimension of the vector space may differ.
  We show that the category of these comodules is, up to a braided abelian equivalence, determined only by the rank of the Hecke symmetry.

  We show that this category is a ribbon category and compute the quantum dimension of simple comodules. Further, we give the explicit formula for the integral on $H_R$.

The work is organized as follows. In Section \ref{sectdef} we briefly recall some definitions of the bialgebra $E_R$ and its Hopf envelope $H_R$. We focuse ourselves on the case of Hecke symmetries. We show that the canonical map $E_R\lora H_R$ in injective, which means that the matrix quantum groups of type $A$ is dense in the corresponding matrix quantum semi-group.
 In Section \ref{sectstruct} we study comodules of $E_R$ and $H_R$ for an even Hecke symmetry $R$. We show that the category $H_R$-\Comod\ is determined by $q$ and the rank of $R$. Further we compute some quantum dimension of simple $H_R$ comodules. In Section \ref{sectint} we compute the integral on $H_R$ and $SH_R$.

In the Appendix we briefly recall some results on Hecke algebra from \cite{dj1,dj2}  and prove some lemmas needed in our context. 

\vskip1ex

\noindent{\it Notation.} Let us fix a field  $\bK$ of characteristic zero. All objects of the paper are defined over $\bK$. Throughout the paper we shall deal with an element $q$ from $\bK^\times$ which we will assume {\it not to be a root of unity, but may be the unity itself.}

 The symmetric groups are denoted by $\Ss_n$, the corresponding Hecke algebras are denoted by $\H_n$. 

A partition $\lambda$ of a non-negative integer $n$ is a sequence of non-increasing non-negative integers whose sum, denoted by $|\lambda|$, is qual to $n$. The length $l(\lambda)$ is the cardinal of its non-zero components. The set of all partitions is denoted by $\P$, the set of partitions to a given $n$ is denoted by $\P_n$. The set of partitions of length $r$ is denoted by $\P^r$. The conjugate partition $\lambda'$ to a given partition  $\lambda$ is the one, whose $i^{\rm th}$ component is $\#\{\lambda_j|\lambda_j\geq i\}$. Some time we shorten the notation of a partition with repeated indices, e.g., $(2,1^3):=(2,1,1,1)$.

Analogously, a $\bZ$-partition is a (finite) sequence of decreasing integers. The set of all $\bZ$-partitions is denoted by $\P^\bZ$.  
\section{The Matrix Quantum Groups and Quantum Semi-groups}\label{sectdef}
\subsection{Yang-Baxter Operators and the Associated Matrix Quantum (Semi-) Groups}
Let $V$ be a finite dimensional vector space over $\bK$ and $R:V\ot V\lora V\ot V$ be an invertible $\bK$-linear operator. Fix a basis $x_1,x_2,\ldots,x_d$ of $V$ and let $\xi^1,\xi^2,\ldots,\xi^d$ be its dual basis in the dual vector space $V^*$ to $V$, $<\xi^i|x_j>=\delta^i_j.$ Denote $z^i_j:=\xi^i\ot x_j$, then $z^i_j,1\leq i,j\leq d$ form a $\bK$-linear basis of $V^*\ot V$. Let $R^{kl}_{ij}$ be the matrix of $R$ with respect to the basis $x_i\ot x_j$: $x_i\ot x_jR=x_k\ot x_lR^{kl}_{ij}$.
The algebra $E_R$ is defined to be the factor algebra of the tensor algebra on $V^*\ot V$ by the relations:
\bba\label{eq0}
R^{mn}_{ij}z^k_mz^l_n=z^p_iz^q_jR^{kl}_{pq},  1\leq i,j,k,l\leq d.\eea
{\it Remark.} The definition of $E_R$ is, in fact, basis free.

$E_R$ is a bialgebra, the coproduct and counit are given by $\Delta(z^i_j)=z^i_k\ot z^k_j$ and $\va(z^i_j)=\delta^i_j$, respectively (we shall frequently use the convention of summing up by the indices that appear both in lower and upper places). We shall be interested in right $E_R$-comodules, the coaction of $E_R$ will be denoted by $\delta$. 

$E_R$ is called {\it matrix quantum semi-group} iff $R$ satisfies the Yang-Baxter equation:
$$ (R\ot \id_V)(\id_V\ot R)(R\ot \id_V)=(\id_V\ot R)(R\ot \id_V)(\id_V\ot R).$$
In this case  $R$ is called {\it Yang-Baxter operator}.

The {\it matrix quantum group} associated to  $R$ is defined to be the {\it Hopf envelope} of $E_R$.
By definition, the Hopf envelope of a bialgebra $E$ is a pair $(H,i)$ consisting of a Hopf algebra $H$ and a bialgebra homomorphism $i:E\lora H$, satisfying the following universal property:
\begin{itemize} \item[] for any Hopf algebra $F$ and a bialgebra homomorphism $e:E\lora F$, there exists  uniquely a Hopf algebra homomorphism $h:H\lora F$, such that $e=h\circ i$. \end{itemize}
This definition is dual to the definition of the maximal subgroup of a semi-group. Because of that we use this definition as the standard way of finding a quantum group from a given quantum semi-group. The Hopf envelope always exists (cf. \cite{manin1,takeuchi}). 

Let $H_R$ be the Hopf envelope of $E_R$, then $H_R$ is, in general, infinitely generated as an algebra. However, with some more condition on $R$, $H_R$ may turn out to be finitely generated. 

Since $R$ satisfies the Yang-Baxter equation, $E_R$ is coquasitriangular, i.e., $E_R$-comod is a braided category (cf. \cite{l-t}). It is natural to ask, whether one can extend the coquasitriangular structure onto $H_R$. 

Let $M$ be an $E_R$-comodule, which is finite dimensional as a vector space. $M$ is naturally endowed with an $H_R$-comodule structure, hence, so is its dual vector space $M^*$. In this case, the braiding on $M\ot M^*$ is uniquely determined by the one on $M\ot M$. In fact, if, with respect to some basis of $M$, the braiding on $M\ot M^*$ has a matrix $T$ then the braiding on $M\ot M^*$ with respect to this basis and its dual should should have a matrix $S$ that satisfies
$$T_{im}^{kn}S_{mj}^{nl}=\delta_i^l\delta_j^k\quad \mbox{and}\quad S_{im}^{kn}T_{nj}^{ml}=\delta_i^l\delta_j^k.$$ 
In particular, this holds for $M=V$ -- the fundamental comodule.
The braiding on $V\ot V^*$ is given by:
\bbs x_l\ot \xi^kR_{V,V^*}=\xi^i\ot x_jP^{kj}_{li},\ees
with $P$ satisfying $R^{il}_{jk}P^{km}_{ln}=\delta^i_n\delta^m_j$ (and hence $P^{il}_{jk}{R^{-1}}^{km}_{ln}=\delta^i_n\delta^m_j$). It was shown by Lyubashenko \cite{lyu1} that in this case, there exists a matrix $Q$ that
 satisfies
${R^{-1}}^{il}_{jk}P^{km}_{ln}=\delta^i_n\delta^m_j$, and the operator twisting $V^*$ and $V$ is given by:
\bbs \xi^i\ot x_jR_{V^*,V}=x_l\ot \xi^k{R^{-1}}^{il}_{jk}.\ees

\begin{thm}\label{thmhr}\cite{frt} If $R$ is a closed Yang-Baxter operator then $H_R$ can be characterized as the factor algebra of the tensor algebra on $z_i^j,t_i^j, 1\leq i,j\leq d$, by the following relations:
\bba R_{ij}^{mn}z^k_mz^l_n&=&z^p_iz^q_jR^{kl}_{pq}\label{eqrzz},\\
t_j^iz^j_k&=&\delta^i_k\label{eqtz},\\
z_j^it^j_k&=&\delta^i_k\label{eqzt}.\eea
The coproduct is: $\Delta(z^j_i)=z^j_k\ot z^k_i, \Delta(t^j_i)=t_i^k\ot t_k^j$,and the antipode is determined by $S(x^i_j)=t^i_j$.\end{thm}
\proof According to the general construction of $H_R$ \cite{manin1},
to show that $H_R$ is a Hopf algebra it is sufficient to show that the antipode on $t^i_j$ is representable in terms of $z_i^j$.

 Multiplying \rref{eqrzz} with $t^s_k$ from the left and $t^j_r$ from the right and using \rref{eqtz}, \rref{eqzt}, we have
\begin{equation} R^{mn}_{ij}z_n^lt^j_q=t_k^mz^p_iR^{kl}_{pq}\label{eqrzt}\end{equation}
or equivalently
\begin{equation}\label{eqpzt}P^{qp}_{lk}z_n^lt^j_q=t_k^mz^p_iP^{ji}_{nm}.\end{equation}
Setting $n=j$ in \rref{eqpzt} and summing up by this index we get
\bba\label{eqtz1}B^p_k=t_k^mz^p_i B^i_m, \eea
where $B^i_m:=P^{li}_{lm}$. Analogously, setting $k=p$ in \rref{eqpzt} and summing up by this index, we get:
\bba\label{eqzt1}z_k^lC^q_lt^j_q=C^j_k,\eea
where $C^i_j:=P^{il}_{jl}$.

The crucial step is to show that $B$ and $C$ are invertible. It turns out that
$$ {C^{-1}}^j_i=Q^{lj}_{li},\quad {B^{-1}}^j_i=Q_{il}^{jl}.$$
The proof of these equalities is based on the Yang-Baxter equation for $R$ and the derived equations of $P$ and $Q$ (cf. Equation \rref{eqinj5}). Hence we have
\begin{equation}\label{eqscd}S(t^j_k)=B^i_kz^p_i{B^{-1}}^j_p={C^{-1}}^n_kz_n^lC^j_l.\end{equation}
Thus we showed that $H_R$ is a Hopf algebra. It this then easy to show that this is the Hopf envelope of $E_R$.\eee

The matrices $B$ and $C$ introduced in the above proof play an important role in the study of comodules of $H_R$. They were called parity operator by Gurevich. In our context, we think the name reflection is more appropriate, since $B$ and $C$ are matrices of the two canonical comodule morphisms $V\lora V^{**}$ and $V^{**}\lora V$.

Equation \rref{eqscd} shows that $V$ being a simple comodule on $E_R$ may turn to be reducible being comodule on $H_R$. In \cite{frt}, there was  introduced a condition ([loc.cit], Eq. (20)), which was equivalent to the following:
\bbs BC=\mbox{ const} \cdot I,\ees
that makes the second equation in \rref{eqscd} trivial.

\subsection{Matrix Quantum (Semi-) Groups Associated with Hecke Symmetries}\label{sectesym}
A Yang-Baxter operator $R$ is called {\it Hecke operator} if it satisfies the following equation
$$(R+1)(R-q)=0, \mbox{ for some } q\in\bK^\times.$$
We shall assume once for all that $q$ is not a root of unity except the unity itself. Further, if  $R$ is closed then $R$ is called {\it Hecke symmetry}.

 We define the  algebras $\Lambd_R$ and $S_R$ as the factor algebras of the tensor algebra upon $V$ by the following relations:
\bba\label{eqslam}\bbar{r} -x_i\ot x_j=x_k\ot x_lR^{kl}_{ij},  1\leq i,j\leq d,\\
qx_i\ot x_j=x_k\ot x_lR^{kl}_{ij}, 1\leq i,j\leq d,\eear\eea
respectively. These algebras play the analogous role of the antisymmetric and symmetric tensor algebras on $V$.

A Hecke symmetry is called {\it even} if the algebra $\Lambda$ defined as above is finite dimensional. The degree of the highest non-vanishing homogeneous component of $\Lambda$ is called the {\it rank} of $R$.

The most well known example of $R$ is the Drinfel'd-Jimbo matrix $R_d$ given by:
 \bba\label{eqmatdj}R_d:=q\sum_{i,j=1}^de_{ij}^{ij}-(q+1)\sum^d_{i,j=1}(\frac{q^{\va(i,j)}}{1+q^{\va(i,j)}}e_{ij}^{ij}-\frac{q^{{\va(i,j)}/2}}{1+q^{\va(i,j)}}e_{ij}^{ji}),\eea
where ${\va(i,j)}:=\mbox{sign} (i-j)$.
$R_{d}$ is an even Hecke symmetry of rank $d=\dim(V)$. The associated Hopf algebra $H_R$ is the general linear quantum group GL$_q(d)$ \cite{frt}. Analogously, the general linear quantum supergroup is given by the following matrix:
\bba\label{eqmatsuper} R_{m|n}=
q\sum^{m+n}_{i,j=1}e_{ij}^{ij}-(q+1)\sum^{m+n}_{i,j=1}(\frac{q^{\va(i,j)}}{1+q^{\va(i,j)}}e_{ij}^{ij}-(-
1)^{\hat{i}\hat{j}}\frac{q^{{\va(i,j)}/2}}{1+q^{\va(i,j)}}e_{ij}^{ji}),\eea
where $\hat{i}$ is the parity of the (homogeneous) basis element $x_i$. $R_{m|n}$ is a non-even Hecke symmetry if $m,n\neq 0$. The associated Hopf algebra $H_R$ is the general linear quantum supergroups GL$_q(m|n)$  \cite{manin2,ls}.

Even Hecke symmetries of rank 2 were classified by Gurevich \cite{gur1}. In particular he showed that on a vector space of dimension at least 2 there exists an even Hecke symmetries of rank 2. Thus, the rank of an even Hecke symmetry may differ from the dimension of the vector space it is defined on.

Another even Hecke symmetry was found by Cremer and Gervais \cite{c-g}.

\subsection{Simple Comodules of $E_R$ and $H_R$}\label{sectcomod}
By its definition, $E_R$ is a graded algebra. Its $n^{\rm{th}}$ component $E_n$ is the factor space of $(V^*\ot V)^{\ot n}$. Therefore, $E_n$ is a coalgebra and every $E_R$-comodule decomposes into the direct sum of $E_n$-comodules \cite{green2}.  $E_n^*$, the dual vector space to $E_n$, is an algebra and acts on every (right) $E_n$-comodules from the left. In fact, we have an equivalence of the two categories: $_{E^*_n}{\cal M}  \cong {\cal M}^{E_n}$.

The Hecke algebra $\H_n$ (see Appendix) acts on $V^{\ot n}$ from the right via $R$, for $x\in V^{\ot n}$,
\bbs xT_i:= x R_i, \mbox{ where } R_i:=\id_{V^\ot i-1}\ot R\ot \id_{V\ot n-i-1}, 1\leq i\leq n-1.\ees
Thus, we have an algebra homomorphism
\bb\label{eqe1}\rho_n:\H_n\lora \End(V^{\ot n}), \rho_n(T_i)=R_i.\end{equation}

In \cite{ph97} the following analogue of Schur's double centralizer theorem was proved.
\begin{thm}\label{schur} The algebras $E^*_n$ and $\rho(\H_n)$ are centralizers of each other in $\End_{\bK}(\Vn)$.\end{thm}

Using standard arguments we deduce form this theorem
\begin{cor} \label{corsimple}Every simple $E_n$-comodule has the form 
\bbs L\cong \Im(\rho(E)),\ees
where $E$ is a primitive idempotent of $\H_n$. Conversely, every primitive idempotent $E$ of $\H_n$ induces a simple comodule of $E_n$ iff $\rho (E)\neq 0$. Further, the modules defined by $E$ and $E'$ are non-isomorphic iff $E$ and $E'$ belong to different minimal two-sided ideals.\end{cor}
 
\noindent{\it Remark.} Since $\H_n$ is the direct product of matrix rings over $\bK$, so is $E_n^*$. Therefore Schur's lemma holds for $E_R$ -- the automorphism ring of a simple  $E_R$-comodule is $\bK$.

\begin{lem} \label{lemsimple} Let $R$ be a Hecke operator. Let $r$ be the largest number, such that $\Lambd_R^r\neq 0$. Assume that $r<\infty$, then $\Lambd_R^r$ is a simple $E_R$-comodule and for every simple $E_R$-comodule $N$, $\Lambd_R^r\ot N$ is simple too.\end{lem}
\proof It is known (cf.  \cite{gur1}) that:
\bbs\Lambd^r_R=\Im(\rho(Y_r)),\ees
where $Y_r$ is the minimal central idempotent that induces the signature representation of $\H_r$, $Y_r=\frac{1}{[r]_{1/q}}\sum_{w\in \Ss_r}(-q)^{-l(w)}T_w$ (see Appendix). Thus, $\Lambd_R^r$ is a simple $E_R$-comodule, since $Y_r$ is a primitive idempotent. By assumption we have:
\bbs \rho(Y_r)\neq 0 \mbox{ and } \rho(Y_s)=0, \forall s\geq r+1.\ees
Let $N$ be a simple $E_R$-comodule of homogeneous degree $m\geq 1$. By Lemma \ref{lemsimple} there exists a primitive idempotent $E$ in $\H_m$, such that $N\cong \Im(\rho(E))$. $\Lambd^r_R\ot N$ is a homogeneous comodule of degree $m+r$. Embedding $\H_r$ and $\H_m$ into $\H_{m+r}$ as follows:
\bbas \H_r\ni T_i\loma T_i\in \H_{m+r}, 1\leq i\le r-1\\
\H_m\ni T_j\loma T_{r+j}\in\H_{m+r},1\leq j\leq m-1,\eeas
we have $\Lambd^r_R\ot N\cong \Im\rho_{m+r}(Y_rE).$ Thus, it is sufficient to show that $\rho_{m+r}(Y_rE)$ is a primitive idempotent in $\rho(H_{m+r})$. This fact will be proved in the Appendix \ref{appprim}.\eee

A direct consequence of this lemma is that $\Lambd_R^r\ot \Lambd_R^r$ is a simple comodule. Hence, the braiding on $\Lambd_R^r\ot \Lambd_R^r$ should be scalar. Since the identity operator is not closed, except when $\dim_\bK\Lambd_R^r=1$, $R$ is not closed if $\dim_\bK\Lambd_R^r\geq 2$ (cf. \cite{gur1}).

Now assume that $R$ is an even Hecke symmetry, then $\dim_\bK\Lambd_R^r=1$, and hence induces the quantum determinant denoted by $D$ in $E_R$. It is shown in \cite{gur1}, that $E_R[D^{-1}]$ is a Hopf algebra, hence coincides with $H_R$ (cf. \cite{haya1}). Let  $\Lambd^{r*}_R$  denote the left dual $H_R$-comodule to $\Lambd^r_R$. Then, by definition of dual comodules and since the antipode sends a group-like element to its inverse,
\bbs \Lambd_R^r\ot\Lambd_R^{r*}\cong\Lambd_R^{r*}\ot\Lambd_R^r\cong I\ees
as $H_R$-comodules.
\begin{edl}\label{dlinj} Let $R$ be a Hecke symmetry. Then the following assertions hold:
\begin{itemize}
\item[(i)] The canonical homomorphism $i:E_R\lora H_R$ is injective,
\item[(ii)] $H_R$ is cosemi-simple,
\item[(iii)] Simple comodules of $H_R$ have the form $N\ot {\Lambd^r_R}^{*\ot n}$, where $N$ is a simple $E_R$-comodule, $n\geq 0$.\end{itemize}\end{edl} 
\proof 
According to Lemma \ref{lemsimple}, $N\ot\Lambd_R^r$ is simple whenever $N$ is. Hence, $D$ is not a zero-divisor in $E_R$. Thus, the localizing homomorphism $E_R\lora E_R[D^{-1}]=H_R$ is injective, (i) is proved.

Let $M$ be a finite dimensional comodule on $H_R\cong E_R[D^{-1}]$, and $x_i,i=1,2,\ldots,m$ be a $\bK$-basis of $M$, then there exist elements $a^i_j,1\leq i,j\leq m$ in $H_R$, such that $\delta(x_i)=x_ja^j_i.$ Let $n$ be such an integer, that $a^i_jD^n\in E_R, \forall i,j$. Then $M\ot {\Lambd^r_R}^{\ot n}$ is an $E_R$-comodule.
 
Let $M\subset N$ be finite dimension $H_R$-comodules. There exists $n\in \Bbb N$ such that $M\ot {\Lambd_R^r}^{\ot n}$ and $N\ot {\Lambd_R^r}^{\ot n}$ are $E_R$-comodules. Hence, there exists $E_R$-comodule $P$, such that:
\bbs N\ot {\Lambd_R^r}^{\ot n}=M\ot {\Lambd_R^r}^{\ot n}\oplus P.\ees
Therefore
\bbs N=M\oplus P\ot {\Lambd_R^r}^{*\ot n}.\ees
Thus every finite dimension $H_R$-comodule is absolutely reducible. Since every $H_R$-comodule is an injective limit of a system of finite dimension comodules, (ii) is proved.
(iii) also follows from this discussion.\eee

\newcommand{\Vect}{\mbox{Vect}}

We now show  that the canonical homomorphism $i:E_R\lora H_R$ is injective for arbitrary Hecke symmetry  $R$.

Let $\V$ be the $\bK$-additive braided category, generated by $V$ and $R$, that is, objects of $\V$ are $V^{\ot n}, n=0,1,2,\ldots,$ and morphisms in $\V$ are obtained from $R$ and the identity morphisms by composing, taking tensor products and linear combinations with coefficient in $\bK$. Analogously, let $\W$ be the $\bK$-additive monoidal category, generated by $V,V^*$ and $R_{V,V},R_{V^*,V},R_{V,V^*}$ and $\ev_V$, $\db_V$. Then $\W$ is a rigid category. In fact, $(V,\ev_V\circ R_{V,V^*})$ is the dual to $V^*$.

We shall use graphical notation (cf.\cite{kt1,kassel}).
The identity morphism $\id_V$ is depicted by an arrow, oriented downward, its dual $\id_{V^*}$ is depicted by an arrow, oriented upward:
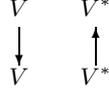
\newcommand{\darr}[2]{\put(#1,#2){\vector(0,-1){0}}}
\newcommand{\uarr}[2]{\put(#1,#2){\vector(0,1){0}}}
\bbs\fbgr{30}{10}\put(0,0){\idgr{10}}\darr{0}{0}
\put(0,0){\ou{V}}\put(0,10){\oo{V}}

\put(20,0){\idgr{10}}\uarr{20}{10}
\put(20,0){\ou{V^*}}\put(20,10){\oo{V^*}}
\egr\ees
The morphism $R_{V,V}$ is depicted by a braid:
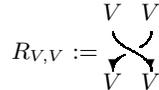
\bbs\fbgr{30}{10}
\put(-5,0){\oo{R_{V,V}:=}}\put(10,0){\ou{V}}
\put(10,10){\oo{V}}\put(20,0){\ou{V}}\put(20,10){\oo{V}}
\put(10,0){\braid} \darr{20}{0}\darr{10}{0}
\egr\ees

To depict composition of morphism we place their pictures in successive order down-ward, the first morphism to be applied will be placed first. To depict tensor product of morphism, we place them in successive order from left to right.

The morphisms $\ev_V$ and $\db_V$ are depicted by special picture, the tangles:
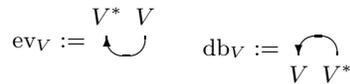
\bbs\fbgr{70}{10}
\put(0,0){\bgr{.25ex}{10}{10}\put(-5,0){\oo{\ev_V:=}}
\put(10,01){\mult}\uarr{10}{06}\put(10,6){\oo{V^*}}\put(20,6){\oo{V}}
\egr}
\put(50,-02){\bgr{.25ex}{10}{10}\put(-5,0){\oo{\db_V:=}}
\put(10,03){\comult}\darr{10}{03}\put(10,1){\ou{V}}\put(20,1){\ou{V^*}}\egr}
\egr\ees

The oriented tangles 
\bbs\mbox{\fbgr{60}{3}
\put(0,-2){\bgr{.25ex}{10}{10}
\put(10,0){\mult}\uarr{20}{05}\egr}
\put(40,-2){\bgr{.25ex}{10}{10}
\put(10,0){\comult}\darr{20}{0}\egr}
\fegr}\ees
are not assumed to present any morphism.

 Reidermeister's moves are allowed, except the second one:
\bbs\fbgr{50}{10}\put(0,0){\mult}\put(0,5){\braid}\put(15,5){\oo{=}}\put(20,0){\mult}\put(20,5){\ibraid}\put(35,5){\oo{=}}
\put(40,0){\mult}\put(40,5){\idgr{10}}\put(50,5){\idgr{10}}\egr\ees
Thus the (oriented) circle does not represent any morphism. Instead, the twisted circles
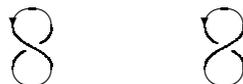
\bbs\fbgr{50}{20}
\put(0,0){\bgr{.25ex}{100}{20}\put(0,0){\mult}\put(0,5){\braid}\put(0,15){\comult}\darr{0}{15}\egr}
\put(50,0){\bgr{.25ex}{100}{20}\put(0,0){\mult}\put(0,5){\ibraid}\put(0,15){\comult}\darr{0}{15}\egr}
\egr\ees
represent (different) morphisms $\bK\lora \bK$.

As a consequence, we obtain the presentation of $R_{V^*,V},R_{V,V^*},R_{V^*,V^*},$

\bbs\mbox{\fbgr{160}{10}

\put(0,0){\bgr{.25ex}{10}{10}
\put(-7,0){\oo{R_{V,V^*}:=}}
\put(10,0){\braid} \darr{20}{0}\uarr{20}{10}
\egr}

\put(60,0){\bgr{.25ex}{10}{10}
\put(-7,0){\oo{R_{V^*,V}:=}}
\put(10,0){\braid} \darr{10}{0}\uarr{10}{10}
\egr}
\put(120,0){\bgr{.25ex}{10}{10}
\put(-10,0){\oo{R_{V^*,V^*}:=}}
\put(10,0){\braid} \uarr{20}{10}\uarr{10}{10}
\egr}
\fegr}\ees

 Since $R$ obeys the Hecke equation, we have

\bba\label{eqinj1}\mbox{\fbgr{160}{30}
\put(0,0){\bgr{.25ex}{10}{10}
\put(0,0){\idgr{15}}\put(10,0){\mult}\put(10,5){\braid}\put(0,15){\braid}\put(0,25){\comult}\put(20,15){\idgr{15}}
\put(0,0){\ou{V^*}}\put(20,30){\oo{V^*}}\uarr{0}{01}\uarr{20}{30}
\put(30,13){\oo{=}}
\egr}
\put(40,-5){\bgr{.25ex}{10}{10}\put(0,05){\idgr{5}}\put(0,10){\ibraid}\put(0,20){\idgr{15}}\put(10,5){\mult}\put(20,10){\idgr{10}}\put(10,20){\braid}\put(10,30){\comult}
\put(0,5){\ou{V^*}}\put(0,35){\oo{V^*}}\uarr{0}{6}\uarr{0}{35}
\put(30,18){\oo{=}}
\egr}
\put(85,0){\bgr{.25ex}{10}{10}\put(0,13){\oo{q^{-1}}}\put(10,0){\idgr{30}}\uarr{10}{30}\put(10,0){\ou{V^*}}\put(10,30){\oo{V^*}}\uarr{10}{30}\egr}
\put(110,0){\bgr{.25ex}{10}{10}\put(03,13){\oo{-(1-q^{-1})}}\put(20,0){\idgr{30}}\put(25,5){\mult}\put(25,10){\braid}
\put(25,20){\comult}\uarr{25}{10}\uarr{20}{30}\put(20,0){\ou{V^*}}\put(20,30){\oo{V^*}}
\egr}
\fegr}\eea

\bba\label{eqinj2}\mbox{\fbgr{160}{30}
\put(0,0){\bgr{.25ex}{10}{10}\put(0,0){\mult}\put(0,5){\braid}\put(0,15){\idgr{15}}\put(10,15){\braid}\put(10,25){\comult}\put(20,0){\idgr{15}}
\darr{0}{29}\darr{20}{0}\put(0,30){\oo{V}}\put(20,0){\ou{V}}\put(30,13){\oo{=}}\egr}

\put(40,0){\bgr{.25ex}{10}{10}\put(0,0){\idgr{15}}\put(0,15){\ibraid}\put(0,25){\idgr{5}}\put(10,0){\mult}\put(10,5){\braid}\put(20,15){\idgr{10}}\put(10,25){\comult}
\darr{0}{29}\darr{0}{0}\put(0,30){\oo{V}}\put(0,0){\ou{V}}\put(30,13){\oo{=}}\egr}
\put(85,0){\bgr{.25ex}{10}{10}\put(0,13){\oo{q^{-1}}}\put(10,0){\idgr{30}}\put(10,0){\ou{V}}\put(10,30){\oo{V}}\darr{10}{0}\egr}
\put(110,0){\bgr{.25ex}{10}{10}\put(03,13){\oo{-(1-q^{-1})}}\put(20,0){\idgr{30}}\put(25,5){\mult}\put(25,10){\braid}
\put(25,20){\comult}\uarr{25}{10}\darr{20}{0}\put(20,0){\ou{V}}\put(20,30){\oo{V}}
\egr}
\fegr}\eea

Further we have:
\bba\label{eqinj3}\mbox{\fbgr{80}{15}
\put(0,3){\bgr{.25ex}{10}{10}
\put(0,0){\braid}\darr{10}{0}\uarr{10}{10}\put(0,0){\ou{V^*}}\put(10,0){\ou{V}}\put(10,10){\oo{V^*}}\put(0,10){\oo{V}}\put(15,3){\oo{=}}\egr}

\put(35,3){\bgr{.25ex}{10}{10}
\put(0,0){\braid}\darr{10}{0}\uarr{10}{10}\put(0,0){\ou{V^*}}\put(10,0){\ou{V}}\put(10,10){\oo{V^*}}\put(0,10){\oo{V}}
\put(-7,3){\oo{q^{-1}}}
\egr}

\put(85,0){\bgr{.25ex}{10}{10}\put(-20,5){\oo{-(1-q^{-1})}}\darr{05}{0}\uarr{30}{15}\put(05,0){\braid}
\put(05,10){\comult}
\put(20,0){\mult}\put(20,5){\ibraid}
\put(05,0){\ou{V^*}}\put(15,0){\ou{V}}\put(30,15){\oo{V^*}}\put(20,15){\oo{V}}
\egr}\fegr}\eea

\bba\label{eqinj4}\mbox{\fbgr{80}{15}
\put(0,2){\bgr{.25ex}{10}{10}
\put(0,0){\braid}\darr{0}{0}\uarr{0}{10}\put(0,0){\ou{V}}\put(10,0){\ou{V^*}}\put(10,10){\oo{V}}\put(0,10){\oo{V^*}}\put(15,3){\oo{=}}\egr}
\put(35,2){\bgr{.25ex}{10}{10}
\put(0,0){\braid}\darr{0}{0}\uarr{0}{10}\put(0,0){\ou{V}}\put(10,0){\ou{V^*}}\put(10,10){\oo{V}}\put(0,10){\oo{V^*}}
\put(-7,3){\oo{q^{-1}}}
\egr}
\put(80,0){\bgr{.25ex}{10}{10}\put(-15,5){\oo{-(1-q^{-1})}}
\put(0,0){\comult}\put(0,10){\mult}\darr{0}{0}\uarr{0}{15}
\put(0,0){\ou{V}}\put(10,0){\ou{V^*}}\put(10,15){\oo{V}}\put(0,15){\oo{V^*}}
\egr}\fegr}
\eea
These equations follow immediately from the equations below:

\bba\label{eqinj5}\mbox{\fbgr{130}{25}
\put(0,0){\bgr{.25ex}{10}{10}\uarr{20}{30}
\put(0,0){\idgr{15}}\put(10,0){\mult}\put(10,5){\ibraid}\put(0,15){\braid}\put(0,25){\comult}\put(20,15){\idgr{15}}
\put(30,13){\oo{=}}
\egr}
\put(40,-5){\bgr{.25ex}{10}{10}\uarr{0}{35}
\put(0,05){\idgr{5}}\put(0,10){\braid}\put(0,20){\idgr{15}}\put(10,5){\mult}\put(20,10){\idgr{10}}\put(10,20){\braid}\put(10,30){\comult}
\put(30,18){\oo{=}}\put(35,5){\idgr{30}}\uarr{35}{35}
\egr}

\put(100,0){\bgr{.25ex}{10}{10}\darr{10}{15}
\put(20,0){\idgr{15}}\put(10,15){\comult}\put(0,10){\mult}\put(0,15){\idgr{15}}\put(25,13){\oo{=}}\put(30,0){\idgr{30}}\uarr{0}{30}\egr}

\fegr}\eea
\begin{edl}\label{dlinja} Let $R$ be a closed Hecke operator then the canonical homomorphism
$i:E_R\lora H_R$ is injective.\end{edl}
\proof
It is known that $E_R$ (resp. $H_R$) is isomorphic to the Coend of the forgetful functor $v:\V\lora\Vect_\bK$ (resp. $w:\W\lora\Vect_\bK)$ (cf. \cite{schauen}). Relations on $E_R$ (resp. $H_R$) are obtained from morphisms in $\V$ (resp. $\W$). In fact, the relations on $E_R$ are obtained from $\Hom_\V(V^{\ot n},V^{\ot m})$. To show that $i:E_R\lora H_R$ is injective, it is sufficient to show that
\bba\label{eqinj0} \Hom_\V(V^{\ot n},V^{\ot m})\cong\Hom_\W(V^{\ot n},V^{\ot m}).\eea
If $m\neq n$, then both sides of \rref{eqinj0} are empty, hence assume that $n=m$.

A tangle which represents a morphism in $\Hom_{\W}(V^{\ot n},V^{\ot n})$ may contain some strands which orientate themselves downward and some closed strands. Using \rref{eqinj3} and \rref{eqinj4} we can disjoint the closed strands, i.e., represent the tangle as a combination of tangle in which the closed strands  are isolated, that is, they do not intersect with other strands. For example, we have:

\bbs\mbox{\fbgr{110}{40}
\put(0,0){\bgr{.25ex}{10}{10}\put(0,0){\idgr{15}}\put(0,15){\braid}\put(0,25){\braid}\put(10,0){\mult}\put(10,5){\braid}\put(20,15){\idgr{20}}\put(0,35){\idgr{5}}\put(10,35){\comult}\darr{0}{39}\darr{0}{25}\uarr{10}{5}\darr{0}{0}\put(25,18){\oo{=}}\egr}
\put(30,05){\bgr{.25ex}{10}{10}
\put(03,13){\oo{q}}\put(10,-5){\idgr{40}}\darr{10}{-5}\darr{10}{34}\put(15,5){\mult}\put(15,10){\braid}\put(15,20){\comult}\uarr{15}{10}
\egr}
\put(83,05){\bgr{.25ex}{10}{10}\put(-13,13){\oo{-(q-1)}}
\put(0,-5){\idgr{20}}\put(10,0){\mult}\put(10,5){\braid}\put(0,15){\braid}\put(0,25){\comult}\put(20,15){\idgr{20}}
\darr{0}{-5}\darr{20}{34}
\egr}

\fegr}\ees

An isolated closed strand represents a morphism of $\bK$, i.e., a scalar. Further, using \rref{eqinj1}, \rref{eqinj2} and \rref{eqinj5} we can transform a non-closed strand into one, that does not have local extremals. This means, a tangle with isolated closed strands is equivalent to a braid up to a scalar. Thus, we see that a tangle which represents a morphism in  $\Hom_{\W}(V^{\ot n},V^{\ot n})$  is a linear combination of braids, hence the following equation holds
\bbs  \Hom_{\W}(V^{\ot n},V^{\ot n})\cong
 \Hom_{\V}(V^{\ot n},V^{\ot n}),\ees 
 from which the theorem follows immediately.\eee




\section{Structure of $H_R$-Comod for an Even Hecke Symmetry $R$}\label{sectstruct}

As it was remarked in the previous section, simple $E_R$-comodules are parameterized by a subset of $\P$ -- the set of all partition. We show in this section that for an even Hecke symmetry $R$ of rank $r$, this set is precisely $\P^r=\{\lambda\in\P|\lambda_{r+1}=0\}$. The tactic is to use the \edim. We show first that an $H_R$-comodule has a non-zero \edim. Further we show that the \edim\ of a simple comodule $M_\lambda$ is non-zero iff $\lambda\in\P^r$.

\subsection{The \edim\ of Simple Comodules}
 Let $\tau$ denote the braiding on $H_R$-Comod. In particular, $\tau_{V,V}=R,\tau_{V^*,V}=R_{V^*,V},\tau_{V,V^*}=R_{V,V^*}$. Let $M$ be an object in $H_R$-Comod, we define its \edim\ to be \cite{majid3}:
\bbs \edim (M):=\ev_M\tau_{M,M^*}\db_M({1_{\bK}}).\ees
For example $\edim (V)=\tr(C).$ A crucial point in our study is the fact, again, due to Gurevich \cite{gur1}, that 
\bba\label{eqqrank} \tr(C)=-[-r]_q.\eea
This quantity is called the quantum rank of $R$.

More general, for a morphism $f:M\lora M$ we define its $\etr$ to be
\bbs \etr(f):=\ev_M\tau_{M,M^*}(f\ot M^*)\db_M({1}_{\bK}).\ees
Then we have 
\bba\label{eqetredim}\edim(M)=\etr(\id_M),\eea
and for morphisms $f:M\lora N$ and $g:N\lora M$,
 \bba\label{eqetr}\etr(fg)=\etr(gf).\eea

Let $M=M_1\oplus M_2$ with morphisms $p_i:M\lora M_i$, $j_i:M_i\lora M$, $p_i\circ j_i={\bf 1}_{M_i}$, $j_1p_1+j_2p_2={\bf 1}_M$. Then 
\bbs\bbar{rl} \edim(M)=&\db_M\tau_{M,M^*}(1\ot(j_1p_1+j_2p_2))\ev_M\\
 =&\tr(j_1p_1)+\tr(j_2p_2)\\
 =&\tr(p_1j_1)+\tr(p_2j_2)\\
 =&\edim(M_1)+\edim(M_2).\eear\ees
Thus the \edim\ is additive.

The idea of using $\edim$ comes from the following lemma.
\begin{lem}\label{lem8tr} Let $M$ be a simple $H_R$-comodule, then $\edim(M)\neq 0.$\end{lem}
\proof Since $M$ is simple,
\bba\label{eqstr1}\Hom(M^*\ot M,\bK)\cong \Hom(M,M)=\bK.\eea
$\db_M:\bK\lora M\ot M^*$ is not trivial, hence there exists a morphism $p:M\ot M^*\lora \bK$, such that $p\circ \db_M={\bf 1}_\bK$. Thus $p\neq 0$, hence $0\neq p\circ\tau^{-1}_{M,M^*}:M^*\ot M\lora \bK$. By \rref{eqstr1} there exists $c\in \bK^\times$, such that $c\cdot p\circ \tau^{-1}_{M,M^*}=\ev_M$, hence $\ev_M\tau\db=c\neq 0$.\eee

Let us fix the following standard rule of defining the dual to the tensor product. Let $V_i$, $i=1,2,\ldots, k$ be $H_R$-comodules, then
\bbs (V_1\ot V_2\ot \cdots \ot V_k)^*:= V^*_k\ot V_{k-1}^*\ot\cdots\ot V^*_1,\ees
with the evaluation map $\ev=\ev_k\ev_{k-1}\cdots\ev_1$. 

Let $M=\Im(\rho(E))$ be a simple $H_R$-comodules, where $E$ is a primitive idempotent in $H_n$.  Then
\bba \label{eqedim-1} \edim(M)=\etr(\rho(E)).
\eea

According to \rref{eqetr}, the map $\rho$ composed with $\etr$ gives us a trace map $\H_n\lora \bK$, denoted by $\tr_r$. We show that this trace map depends actually only on $q$ and $r=\rank R$ and not on the operator $R$ itself.

For a morphism $f\in\End^{E_R}(\Vn)$ we define a morphism $\etr^n(f)\in\End^{E_R}(V^{\ot n-1})$ to be:
\bba\label{eq-tr0}\etr^n(f):=(\id_{V^{\ot n-1}}\ot \ev_v\circ\tau_{V,V^*})(f\ot\id_V)(\id_{V^{\ot n-1}}\ot\db).\eea
Thus, we have a map $\End^{E_R}(\Vn)\lora\End^{E_R}(V^{\ot n-1})$.

\begin{lem}\label{lem100} Let $f\in\End^{E_R}(\Vn)$. Then
\bba \etr(f)=\etr^1\circ\etr^2\circ\cdots\circ\etr^n(\tau_{w_n}^{-2}\circ f),\eea
where $w_n$ is the longest element of $\Ss_n$ -- the one, that reverses the order of the elements $1,2,\ldots, n$. \end{lem}
\proof
The above equation above can be best explained using pictorial representation. Below is the picture for $n=3$.
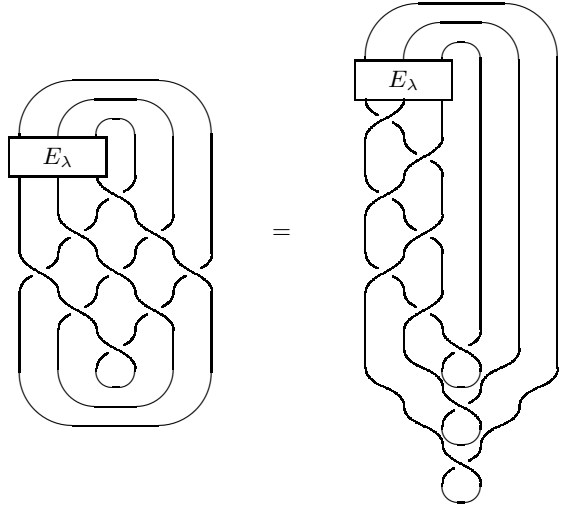
\bbs
\fbgr{160}{130}
\put(0,20){
\bgr{.25ex}{100}{100}
\put(25,15){\oval(50,30)[b]}\put(25,15){\oval(30,20)[b]}\put(25,15){\oval(10,10)[b]}
\put(0,15){\idgr{20}}\put(0,35){\braid}\put(0,45){\idgr{20}}\put(0,65){\multmor{E_\lambda}{20}{10}}\put(10,15){\idgr{10}}\put(10,25){\braid}\put(10,45){\braid}\put(10,55){\idgr{10}}\put(20,15){\braid}\put(20,35){\braid}\put(20,55){\braid}\put(30,25){\braid}\put(30,45){\braid}\put(30,65){\idgr{10}}\put(40,15){\idgr{10}}\put(40,35){\braid}\put(40,55){\idgr{20}}\put(50,15){\idgr{20}}\put(50,45){\idgr{30}}
\put(25,75){\oval(50,30)[t]}\put(25,75){\oval(30,20)[t]}\put(25,75){\oval(10,10)[t]}
\egr}\put(70,70){=}

\put(90,30){
\bgr{.25ex}{150}{100}
\put(25,85){\oval(50,30)[t]}\put(25,85){\oval(30,20)[t]}\put(25,85){\oval(10,10)[t]}
\put(0,15){\idgr{10}}\put(0,25){\ibraid}\put(0,35){\idgr{10}}\put(0,45){\ibraid}\put(0,55){\idgr{10}}\put(0,65){\ibraid}\put(0,75){\multmor{E_\lambda}{20}{10}}\put(10,15){\ibraid}\put(10,35){\ibraid}\put(10,55){\ibraid}\put(20,25){\idgr{10}}\put(20,45){\idgr{10}}\put(20,65){\idgr{10}}\put(30,15){\idgr{70}}\put(40,15){\idgr{70}}\put(50,15){\idgr{70}}
\put(20,5){\braid}\put(20,0){\mult}
\put(20,-0){\twist{-10}{10}}\put(10,10){\idgr{5}}\put(30,0){\twist{10}{10}}\put(40,10){\idgr{5}}\put(20,-10){\braid}\put(20,-15){\mult}
\put(20,-15){\twist{-10}{10}}\put(10,-5){\twist{-10}{10}}\put(0,5){\idgr{10}}
\put(30,-15){\twist{10}{10}}\put(40,-5){\twist{10}{10}}\put(50,5){\idgr{10}}

\put(20,-25){\braid}\put(20,-30){\mult}
\egr}\egr\eqeee
\ees

Let $w\in\Ss_n$ be expressed in the form $w=v_kv_{k+1}\cdots v_{n-1}w_1$, $w_1\in\Ss_{n-1}$. Then it is easy to check
\bba\label{eq-tr1}
\etr^n(\rho(T_w))=\left[\bbar{ll}
\rho_{n-1}(T_{v_k\cdots v_{n-2}w_1})&\mbox{if } k\leq n-1\\
\tr(C)\rho_{n-1}(T_{w_1})&\mbox{if }k=n.\eear\right.\eea
Thus we see that $\etr^n$ depends only on $\tr(C)=-[-r]_q.$ Let us define an analogous operator $\tr^n:\H_n\lora\H_{n-1}$, setting
\bba\label{eq-tr2}
\tr_r^n(T_w)=\left[\bbar{ll}
 T_{v_k\cdots v_{n-2}w_1}&\mbox{if } k\leq n-1\\
-[-r]_qT_{w_1}&\mbox{if }k=n,\eear\right.\eea
where $w=v_kv_{k+1}\cdots v_{n-1}w_1, w_1\in\Ss_{n-1}.$
Then we have
\bba\label{eqtretr}\rho_{n-1}\circ\tr^n_r=\etr^n\circ\rho_n.\eea
Recall that the trace map $\tr_r$ is defined as the composition of $\rho_n$ by $\etr$.
\begin{cor}\label{cortr} We have
\bba\label{eqtr3}\tr_r(E)=\tr_r^1\circ\tr_r^2\circ\cdots\circ\tr_r^n(T_{w_n}^{-2}E).\eea
Hence $\tr_r$ depends only on $r$.\end{cor}
\proof For $n=0$, we have $\rho_0:\H_0\cong\bK\cong\End^{E_R}(V^{\ot 0}).$
Using equation \rref{eqtretr}, we have
\bbas \tr_r(W)&=&\etr(\rho_n(W))\\
&=& \etr^1\circ\etr^2\circ\cdots\circ\etr^n(\tau^{-2}_{w_n}\rho_n(W))\\
&=&\rho_0(\tr_r^1\circ\tr_r^2\circ\cdots\circ\tr_r^n(T^{-2}_{w_n}W))\\
&=&\tr_r^1\circ\tr_r^2\circ\cdots\circ\tr_r^n(T^{-2}_{w_n}W).\eqeee\eeas

\begin{pro}\label{pror+1} Let $R$ be an even Hecke symmetry of rank $r$. Let $\lambda\part n$ and $M=\Im\rho_n(E_\lambda)$. Then $M_\lambda\neq 0$ (and hence is a simple comodule of $H_R$) iff $\lambda_{r+1}=0$.\end{pro}
\proof 
By virtue of Lemma \ref{lem8tr}, it is sufficient to show that
\bbs \tr_r(T^{-2}_{w_n}E):=\tr^1_{r}(\tr^2_{r}(\cdots(\tr^n_{r}(E))\cdots))\neq 0\quad\iff\quad\lambda_{r+1}=0.\ees
 Let $R^{DJ}$ be the Drinfel'd-Jimbo $R$-matrix of type $A_{r-1}$, defined in \rref{eqmatdj}, then $R^{DJ}$ has rank $r$. Let $\rho^{DJ}_n$ denote the representation of $\H_n$ induced by $R^{DJ}$. If for a partition $\lambda$, $\lambda_{r+1}=0$, then $V^{\ot n}$, being considered as $\H_n$-module, contains direct summand isomorphic to the simple $\H_n$-module $S_\lambda$ (see Appendix) (\cite{dpw}, Proposition 5.1). On this module $E_\lambda\neq 0$, hence $\Im(\rho^{DJ}_n(E_\lambda))\neq 0$. By Lemma \ref{lem8tr}, $\tr_{r}(T^{-2}_{w_n}E_\lambda)\neq 0$.

If $\lambda_{r+1}\neq 0$, then $V^{\ot n}$, being considered as $\H_n$-module through $\rho_n^{DJ}$, does not contain direct summand that is isomorphic to $S_\lambda$ as $\H_n$-module, so that $E_\lambda$ acts as zero on $V^{\ot n}$. Hence $ \tr_{r}(T^{-2}_{w_n}E_\lambda)$ acts as zero on $V^{\ot 0}=\bK$. On the other hand, $ \tr_{r}(T^{-2}_{w_n}E_\lambda)$ is a constant, hence it is equal to zero.\eee
\subsection{The Grothendieck Ring of $H_R$-Comod}
Let $E_\lambda$ and $E_\mu$ be primitive idempotents of $\H_n$ and $\H_m$ respectively, $M_\lambda=\rho_n(E_\lambda), M_\mu=\rho_m(E_\mu)$. Embed $\H_n$ into $\H_{n+m}$, $T_i\loma T_i$ and $\H_n$ into $\H_{n+m}$, $T_j\loma T_{j+n}$. Then $ M_\lambda \ot M_\mu\cong \Im(\rho_{n+m}(E_\lambda E^n_\mu))$, where $E^n_\mu$ is the image of $E_\mu$ under the defined above embedding. The decomposition of $E_\lambda E^n_\mu$ into primitive idempotents depends only on $\lambda$ and $\mu$ themselves.  Consequently, in the decomposition
\bba\label{decomp}M_\lambda\ot M_\mu=\bigoplus_{\gamma\part|\lambda|+|\mu|}c_{\lambda \mu}^\gamma M_\gamma\eea
the coefficients $c_{\lambda \mu}^\gamma$ coincide with the one for (rational) representations of GL$(r)$, i.e., the Littlewood-Richardson coefficients.

Hence we have the determinantal form for $M_\lambda$ (cf. \cite{james1}, Section 19 ):
\bba\label{eqdetform}M_\lambda\cong \det\left| \left(M_{(\lambda_i-i+j)}\right)_{1\leq i,j\leq r}\right|,\eea
where, in the right-hand side, the determinant should be understood to be taken in the Grothendieck ring of $H_R$-Comod and $M_{(k)}:=0$ whenever $k<0.$

We now characterize the dual comodule to $M_\lambda$, thus furnish all $H_R$-comodules. Call sequence of $r$ integers $\lambda=(\lambda_1,\lambda_2,\ldots,\lambda_r), \lambda_i\in\bZ,\lambda_1\geq \lambda_2\geq\ldots\geq\lambda_r$ a $\BBb{Z}$-partition. Let  $\ove{\lambda}_i:=\lambda_i-\lambda_r$. Define
\bbas M_\lambda:=M_{\ove{\lambda}}\ot \Lambda^{r\ot \lambda_r},\eeas
where $\ove{\lambda}=(\ove{\lambda}_1,\ove{\lambda}_2,\ldots,\ove{\lambda}_r)$, and 
\bbas \Lambda^{r\ot \lambda_r}:=\left[\bbar{cc}\Lambda^{r\ot \lambda_r}& \mbox{ if } \lambda_r\geq 0;\\[1ex]\Lambda^{r*\ot |\lambda_r|}& \mbox{ otherwise. }\eear\right.\eeas

For a  $\bZ$-partition $\lambda$, define $-\lambda:=(-\lambda_r,-\lambda_{r-1},\ldots,-\lambda_1)$. We show that $M_\lambda^*\cong M_{-\lambda}$. In fact, since $(M^*_\lambda)^*\cong M_\lambda$, we have
\bbas \Hom(M^*_\lambda,M_{-\lambda})&\cong& \Hom(\BBb{K}, M_{-\lambda}\ot M_\lambda)\\
&\cong& \Hom(\Lambda^{r\ot \lambda_1},M_{\ove{\lambda}}\ot M_\lambda)\\
&\cong&  M(c_{\ove{\lambda},\lambda}^{(\lambda_1^r)})
\cong \bK,\eeas
where $M(c_{\ove{\lambda},\lambda}^{(\lambda_1^r)})$ denotes the matrix ring of rank $c_{\ove{\lambda},\lambda}^{(\lambda_1^r)}$. The last equation above holds, since it does not depend on $R$ and it holds in the classical case of GL$(n)$. One can also use Littlewood-Richardson's rule (cf. \cite{mcd2}) to show directly that  $c_{\ove{\lambda},\lambda}^{(\lambda_1^r)}=1$.

\begin{edl}\label{dlfusion}Let $R$ be an even Hecke symmetry of rank $r$. Then simple $E_R$ comodules are parameterized by partitions of length $\leq r$ and simple $H_R$-comodules are parameterized by $\Bbb Z$-partitions of length $\leq r$. 
 For $\lambda\part n$, $\mu\part m$, we have
\bba\label{eqfusion}M_\lambda\ot M_\mu\cong\bigoplus_{\gamma\in\P^{n+m}_r}c^\gamma_{\lambda,\mu}M_\gamma,\eea
where $c_{\lambda,\mu}^\gamma$ are the Littlewood-Richardson coefficients. The dual to $M_\lambda$ comodule is $M_{-\lambda}$.
Hence, the category ($E_R$-) $H_R$-Comod is equivalent as a braided abelian category to the category of representations of the standard quantum (semi-) group (M$(r)$) GL$_q(r)$.
\end{edl}
\proof It remains to prove the last statement. Remark that every $E_R$-comodule is a submodule of $\Vn$, $n=1,2,\ldots$ and $\End^{E_R}(\Vn)$ is the factor algebra of $\H_n$ by the ideal generalized by minimal central idempotents that correspond to partitions of length larger then $r$. Thus, the category $E_R$-comod, up to braided equivalence, does not depend on $R$ but only on its rank $r$. 

Analogous discussion holds for $H_R$-comodules.\eee

From the above theorem one derives immediately an analog of Peter-Weyl's decomposition for $H_R$. We have an isomorphism
\bba\label{eq-pw}\bigoplus_{\lambda\in\P^\bZ\atop l(\lambda)\leq r}\End_\bK(M_\lambda)\cong H_R.\eea

\subsection{The Computation of \edim\ and \rdim\ of Simple Comodules}
The \edim\ is an invariant in $H_R$-Comod, however it is not tensor multiplicative and hence is difficult to compute for non-simple comodules. Even the formula for $\edim(V^{\ot n})$ is very complicated. There is a so called ribbon-dimension -- rdim -- which is multiplicative and, on simple comodules, it differs from \edim\  only by a power of $q$. In fact, for a simple $E_R$-comodule $M$ of degree $n$, defined by a primitive idempotent $E\in \H_n$, its ribbon-dimension is, by definition,
\bba\label{eqrdim1}\rdim(M):=q^{n(r+1)}\tr_r(T^2_{w_n}E)=q^{n(r+1)}\tr^1\circ\tr^2\circ\cdots\circ\tr^n(E).\eea

By definition, a ribbon in a strict braided category is a natural isomorphism $\theta$, subject to the following equations:
\bba\label{eqrdim2}\bbar{c}\theta_{M\ot N}=(\theta_M\ot\theta_N)\tau_{N,M}\tau_{M,N},\\
\theta_{M^*}=(\theta_M)^*,\eear\eea
where $\tau$ denotes the braiding in the category. This definition first appeared in the works of Joyal and Street \cite{js3} and Turaev \cite{tur1}. A strict braided rigid category with a ribbon is called ribbon category. For more detail on ribbon categories, the reader is referred to a book of C. Kassel \cite{kassel}. For an object $M$ in a ribbon category, the ribbon-dimension, $\rdim(M)$, is by definition the morphism
\bba\label{eqrdim3}\ev_{M}\tau_{M,M^*}(\theta_M\ot M^*)\db_M:I\lora I.\eea
By virtue of \rref{eqrdim2}, \rdim\ is multiplicative (cf. \cite{kassel}, Chapter XIV). 

\begin{lem}\label{lemrdim1} $H_R${\rm -Comod} possesses a ribbon $\theta$, such that for an   idempotent  $E$ in $\H_n$, and $M=\Im(\rho(E))$, the ribbon on $M$ is given by
\bba\label{eqrdim4}\theta_M=q^{n(r+1)/2}{T_{w_n}}^{2}E.\eea
Thus, $H_R$-{\rm Comod} is a ribbon category.\end{lem}
\proof 
Assume that $\theta$ is a ribbon in $H_R$-\Comod, then (cf. \cite{rt})
\bbs \theta^{-2}_M=(\ev_M\tau_{M,M^*}\ot\id_M)(\id_M\ot \tau_{M,M^*}\db_M),\ees
for any comodule $M$. Thus, using Equation \rref{eqinj1}, we have $\theta_V=q^{(r+1)/2}\id_V$. Hence, according to the first axiom for $\theta$, for an idempotent $E\in\H_n$ and $M=\Im\rho(E)\subset\Vn,$
\bbs \theta_M=q^{n(r+1)/2}\rho(T^2_{w_n}E).\ees
Remark that ${T_{w_n}}^{2}$ is central in $\H_n$. Further, we have 
$\theta_{\Lambda_n}=q^{\frac{n(r+1)}{2}}\id_{\Lambda^n}.$
Hence we set 
$\theta_{\Lambda^{r *}}=q^{\frac{r(r+1)}{2}}\id_{\Lambda^{r* }}$
and extend it on the whole category using Theorem \ref{dlinj}, (iii) and the first axiom in \rref{eqrdim2}. 

It remains to verify $\theta_{M^*}=(\theta_M)^*$. We can assume $M\subset\Vn$ hence reduce it to checking $\theta_{\Vn^*}=(\theta_\Vn)^*$, and then to checking $\theta_{V^*}=(\theta_V)^*$, by using the identity
\bba\label{identity}\tau_{V\ot W,U}\circ\tau_{U,V\ot W}=(\tau_{U,V}\ot W)(V\ot(\tau_{W,U}\circ\tau_{U,W}))(\tau_{U,V}\ot W).\eea
We have $V^*\cong\Lambda^{r-1}\ot \Lambda^{r*}$. The problem reduces to checking
$$(\theta_{\Lambda^{r-1}}\ot \theta_{\Lambda^{r*}})\circ \tau_{\Lambda^{r*},\Lambda^{r-1}}\circ\tau_{\Lambda^{r-1},\Lambda^{r*}}=\id_{V^*}.$$
Or, equivalently
$$\tau_{\Lambda^r,\Lambda^{r-1}}\circ\tau_{\Lambda^{r-1},\Lambda^r}=q^{r(r+1)}\id_{\Lambda^{r-1}\ot\Lambda^r}^*.$$
Again, this equation follows from
\bbas \tau_{V^{\ot n},\Lambda^{r}}\tau_{\Lambda^{r},V^{\ot n}}=q^{(r+1)n}\id_{\Lambda^{r},V^{\ot n}}.\eeas
By means of Identity \rref{identity}, an induction reduces it to the check of the above equation for the case $n=1$. We have
\bba\tau_{V^{\ot n},\Vn}^2&=&
R_rR_{r-1}\cdots R_1R_1R_2\cdots R_r\nonumber\\
&=&q^r+(q-1)(q^{r-1}R_r+q^{r-2}R_{r-1}R_rR_{r-1}+\cdots+R_1\cdots R_r\cdots R_1).\label{eqtau2}\eea
Hence
$$\tau_{V,\Lambda^{r}}\tau_{\Lambda^{r},V}=(\rho(Y_r)\ot\id_V)\tau_{V^{\ot n},\Vn}^2=q^{r+1}\id_{\Lambda^r\ot V}.$$
Here we use the fact that $\rho(Y_{r+1})=0$.
\eee

Using the definition in \rref{eqrdim3}, we can now define, for $M=\rho(E)$, where $E$ is an idempotent,
\bba\label{eqrdim41}
\rdim(M)=q^{n(r+1)/2}\tr_r(T_{w_n}^2E)=q^{n(r+1)/2}\tr^1\circ\tr^2\circ\cdots\circ\tr^n(E).\eea

According to \rref{eqdetform} and since \rdim\ is tensor multiplicative, we have the formula for the \rdim\ of $M_\lambda$:
\bba\label{eqrdim5}\rdim(M_\lambda)=\det\left(\rdim(M_{(\lambda_i-i+j)})\right)_{1\leq i,j\leq r}.\eea
Direct computation shows that
\bba\label{eqrdim6}\rdim(M_{(k)})=q^{k(r+1)/2}\tr_{r}^{1,2,\ldots,k}(E_{(k)})=q^{-k(r-1)/2}\frac{[k+r-1]_q}{[r-1]_q![k]_q!}.\eea
Using Ex. I.3.1 of \cite{mcd2} we deduce
\begin{edl}\label{prodim1}Let $M_\lambda$ be a simple $E_R$-comodule corresponding to a partition $\lambda$. Then
\bba\label{eqstar2}\rdim(M_\lambda)=q^{|\lambda|(r+1)/2-\sum_{i=1}^r\lambda_i(r+1-i)}\prod_{1\leq i<j\leq r}\frac{[\lambda_i-\lambda_j+j-i]_q}{[j-i]_q}.\eea
\end{edl}
Since $\rdim(M_{(1^r)})=1$ and since the $\rdim$ is tensor multiplicative, Formula \rref{eqstar2} holds for any $\bZ$-partition. As a direct corollary of this theorem we have
\begin{pro}\label{dldhc} 
Let $\lambda$ be a $\BBb{Z}$-partition. Then
\bba\label{eqdhc} \edim(M_\lambda)=q^{-\frac{1}{2}(|\lambda|^2+\sum_{i=1}^r\lambda_i(\lambda_i+2r-4i+2))}
\prod_{1\leq i<j\leq r}\frac{[\lambda_i-\lambda_j+j-i]_q}{[j-i]_q}.\eea
\end{pro}
\proof It is sufficient to show that, for $n=|\lambda|$,
\bba\label{eqedim0} T_{w_n}^{-2}E_\lambda=q^{-(\sum_i\lambda_i(\lambda_i-2i+1)+n(n-1))/2}E_\lambda.\eea
 This is the assertion of Corollary \ref{cortwn} in Appendix.\eee
\subsection{The Special Matrix Quantum Groups $SH_R$} 
Assume now that the quantum determinant commutes with elements of $H_R$ (see \cite{gur1} for a necessary and  sufficient condition). Then we can set it equal to 1 to get the special matrix quantum  group $SH_R$. 

\begin{thm}\label{thmshr} The Hopf algebra $SH_R$ is co-semi-simple and its simple comodule are parameterized by partitions of length $\leq \rank(R)-1$.
\end{thm}
\proof In the isomorphism \rref{eq-pw}, assume that $H_\lambda$ is the image of $\End_\bK(M_\lambda)$. Then $H_\lambda$ is a simple subcoalgebra of $H_R$. The multiplication on $H_R$ maps
$$ H_\lambda\ot H_\mu\lora \bigoplus_{\gamma\part|\lambda|+|\mu|\atop l(\gamma)\leq r}H_\gamma,\quad (r=\rank(R).)$$
In particular, it maps $H_\lambda\ot H_{(1^r)}\lora H_{\lambda\cup(1^r)}$, where $\lambda\cup(1^r)$ is the partition $(\lambda_1+1,\lambda_2+1,\ldots,\lambda_r+1)$, which is an isomorphism.

Recall that $H_{(1^r)}$ is spanned by the quantum determinant. Thus, if we set the quantum determinant equal to 1, \rref{eq-pw} reduces to the isomorphism
\bba\label{eq-spw}\bigoplus_{l(\lambda)\leq r-1}H_\lambda\cong SH_R.\eea
To see this, we consider the composition
$$\bigoplus_{l(\lambda)\in P^\bZ\atop l(\lambda)\leq r}H_\lambda\lora H_R
\lora SH_R$$ which is obviously surjective. For a partition $\lambda$, let $\overline{\lambda}=(\lambda_1-\lambda_r,\lambda_2-\lambda_r,\ldots,0)$. The discussion above shows that $H_\lambda$ and $H_{\overline{\lambda}}$ have the same image in $SH_R$. Further, since the determinant is not a zero-divisor and lies in the center of $H_R$, the above map in injective when restricted on each $H_\lambda$. Whence \rref{eq-spw} is an isomorphism. \eee

The comodule $\Lambda^r$, being considered as a $SH_R$-comodule, is isomorphic to the trivial module. Hence, it is expected that this module has an \edim\ equal to 1. This can be made by rescaling $R\loma q^{-\frac{r+1}{2r}}R$. We call the rescaled \edim\ the normalized \edim\ and denote it by $\ove{\edim}$,
\bbas \ove{\edim}(M_\lambda)= q^{-\frac{1}{2}(\frac{n^2}{r}+\sum_{i=1}^r\lambda_i(\lambda_i+2r-4i+2))}\prod_{1\leq i<j\leq r}\frac{[\lambda_i-\lambda_j+j-i]_q}{[j-i]_q}\\
=q^{-\frac{1}{2r}\sum_{1\leq i<,j\leq r}(\lambda_i-\lambda_j)^2-
\sum_{i=1}^r\lambda_i(r-2i+1)}
\prod_{1\leq i<j\leq r}\frac{[\lambda_i-\lambda_j+j-i]_q}{[j-i]_q}\eeas
The last equation shows particularly that the $\ove{\edim}$ of $M_\lambda$ and of $M_\lambda\ot \Lambda^r$ are equal. As a corollary of the above discussion we have
\begin{pro}\label{prosl} Let $R$ be such an even Hecke symmetry that the quantum determinant commutes with all elements of $H_R$, so that one can set it equal to 1 to get the corresponding special matrix quantum groups $SH_R$. Then the coquasitriangular structure $<\cdot|\cdot >$ on $SH_R$ can be defined, setting
\bbs <z_i^k|z_j^l>=q^{-\frac{r+1}{2r}}R_{ij}^{lk}.\ees\end{pro}

The ribbon in $SH_R$-\Comod\ is then given by
\bbs\ove{\theta}_{\Vn}=q^{\frac{r(r^2-1)}{2r}}\ove{\tau}_{w_n}^2,\ees
$\ove{\tau}:=q^{-\frac{r+1}{2r}}R.$ Direct computation shows that the  \rdim\ remains unchanged.

From all results obtained we can interpret the quantum group of type $A_n$ as follows.
\begin{dn}\label{qgan}\rm The quantum group ${\bf GL}_{q}({r})$ (${\bf SL}_q(r)$) is a $\bK$-linear braided rigid semi-simple abelian category, whose simple objects are parameterized by $\Bbb Z$-partitions $\{\lambda|\lambda_{r+1}=0\}$ (partitions $\{\lambda|\lambda_r=0\}$) and the tensor product satisfy Equation \rref{eqfusion}. The braiding on the simple object $V$ that corresponds to partition $(1)$ satisfies the equation $(R-q)(R+1)=0$.
\end{dn}
Each even Hecke symmetry of rank $r$ induces a (monoidal) functor from ${\bf GL}_{q}({r})$ into Vect$_\bK$, called a realization of this quantum group. For a realization of ${\bf SL}_q(r)$ we need certain extra condition on the Hecke symmetry.

\noindent{\it Remark.} According to a result of Kazhdan-Wenzl \cite{k-w}, each monoidal category with the Grothendieck ring isomorphic to the one of the category of rational representations of GL$(r)$ (SL$(r)$) should be equivalent to ${\bf GL}_{q}({r})$ (${\bf SL}_q(r)$) for some $q$.

According to Theorem \ref{prodim1}, we have
\begin{pro}For $q\neq q'$, ${\bf GL}_{q}({r})$ and ${\bf GL}_{q'}({r})$ (resp. ${\bf SL}_q(r)$ and ${\bf SL}_{q'}(r)$) are not equivalent as braided categories.\end{pro}
Indeed, if there were such an equivalent, the \rdim\ of the corresponding objects would be equal, which is obviously not the case if $q\not= q'$.

\section{The Integral}\label{sectint}
 By definition, a (left) integral on a Hopf algebra $H$ is a linear functional $\int:H\lora k$ satisfying the following condition:
\bb\label{eqi1}\mbox{ if }\Delta(x)=\sum_{(x)}x_1\ot x_2\quad\mbox{then}\quad \int(x)=\sum_{(x)}x_1\int(x_2).\end{equation}
Since $H_R$ and $SH_R$ are co-semi-simple, such an integral on them exists. If we rescale it, setting $\int(1)=1$, then it exists uniquely.
 Our aim in this section is to compute explicitly the integral for $H_R$ and $SH_R$.
\subsection{The Integral on $H_R$}
Let us make $H_R$ into a $\Bbb Z$-graded algebra by setting $\mbox{deg}_\bZ~(z_i^j)=1,\mbox{deg}_\bZ~(t_i^j)=-1 $. For $x$ homogeneous, we can choose an expression $\Delta(x)=x_1\ot x_2$, such that $\mbox{ deg}_\bZ(x)=\mbox{ deg}_\bZ(x_1)=\mbox{ deg}_\bZ(x_2)$. Hence, \rref{eqi1} is not a  homogeneous relation unless $\mbox{deg}_\bZ(x)=0$, while $H_R$ is homogeneous. Thus, we proved the following.
\begin{lem}\label{lemhomo} The value of $\int$ on a homogeneous element of $H_R$ (with respect to the above grading) is equal to zero unless the degree of this element is zero.\end{lem}

The problem reduces to finding the values of $\int$ on monomials in $z^i_j$'s and $t^k_l$'s, in which the numbers of $z^i_j$'s and $t^k_l$'s are equal. By virtue of the relation in \rref{eqrzt}, every monomial in $z^i_j$'s and $t^k_l$'s can be expressed as a linear combination of monomials in which $z^i_j$'s are on the right of $t^k_l$'s. Hence it reduces to computing:
\bbs \int(\mono{z}{i}{j}\mono{t}{k}{l}).\ees
At first sight, it is not clear if a set of values of $\int$ on $\mono{z}{i}{j}\mono{t}{k}{l}$ induces an integral on $H_R$, since there may be more than one way of expressing a monomial in $z$'s and $t$'s as a linear combination of $\mono{z}{i}{j}\mono{t}{k}{l}$. Here the Yang-Baxter plays an important role, which ensures that different ways of expressing yield the same result. 

We need some notations.
The Murphy operators $L_m$, $1\leq m\leq n$, were introduced in \cite{dj2} as follows: $L_1:=0\footnote{There was a misprint in \cite{dj2} where $L_1:=1$.},$
\bbas 
 L_m:= q^{-1}T_{m-1}+q^{-2}T_{m-2}T_{m-1}T_{m-2}+\cdots+q^{-m+1}T_1T_2\cdots T_{m-1}\cdots T_1,m\geq 2.\eeas

For an element $W$ of $\H_n$, we shall denote $\ove{W}$ it image under the map $\rho:\H_n\lora \End^{E_R}(\Vn)$, the image of $\H_n$ itself will be denoted by $\ove{\H_n}$. Lemma \ref{lemtrq} in the Appendix shows that the operator $(\Ln{n}-[-r]_q)$ is invertible in $\ove{\H_n}=\End^{E_R}(\Vn).$
\begin{edl} \label{dlhs}  Let $R$ be an even Hecke symmetry of rank $r$. Then the integral of $H_R$ can be given by the following formula 
\bba\label{eqi16} \int Z^J_IT^L_K=\sum_{w\in \Ss_n}q^{-l(w)} \left(\prod_{k=1}^n(\Ln{k}-[-r]_q)^{-1}R_{w^{-1}}C^{\ot n}\right)_I^{L'} {R_w}_{K'}^J.\eea
where $I,J,K,L$ are multi-indices, $Z_I^J:=\mono{z}{i}{j}$, and for $K=(k_1,k_2,...,k_n)$, $K':=(k_n,k_{n-1},...,k_1)$. The matrix $C$ was introduced in Equation \rref{eqzt1}.\end{edl}

\proof We should check that $\int$ satisfies \rref{eqi1} and is well defined, i.e., is compatible with the defining relations. The verification of \rref{eqi1} is rather trivial. Let us check that the formula for $\int$ is compatible with the relations \rref{eqtz} and \rref{eqzt}. For simplicity let us denote $P_n:=\prod_{k=1}^n(L_{k}-[-r]_q)$. Thus $\ove{P}_n=\rho_n(P_n)$ is invertible in $\ove{\H_n}$.

Remark that the equation \rref{eqtz} together with \rref{eqpzt} implies \rref{eqzt1}. Thus we have to check
\bba \sum_{w\in\Ss_n}q^{-l(w)}(\overline{P_n}^{-1}R_{w^{-1}}C^{\ot n})_{I_1m}^{L_1m} {R_w}_{K}^J=
\sum_{w\in\Ss_{n-1}}q^{-l(w)}(\ove{P_{n-1}}^{-1}R_{w^{-1}}C^{\ot n-1})_{I_1}^{L_1}{R_w}_{K_1}^{J_1}\delta_{k_n}^{j_n},\label{eq22}\eea
\bba
\sum_{w\in\Ss_n}q^{-l(w)}(\ove{P_n^{-1}}R_{w^{-1}}C^{\ot n-1})_{I_1i_n}^{L_1l_n} {R_w}_{K_1p}^{J_1n}C_n^p=\sum_{w\in\Ss_{n-1}}q^{-l(w)}(\ove{P_{n-1}}^{-1}R_{w^{-1}}C^{\ot n-1})_{I_1}^{L_1}{R_w}_{K_1}^{J_1}C_{i_n}^{l_n},\label{eq23}\eea
here, for a multi-index $K=k_1,k_2,\ldots,k_n$ we denote $K_1:=k_1,k_2,\ldots,k_{n-1}$.

First, we need a technical lemma.
\begin{lem}\label{technicallem} The following relation holds in $\H_n$:
\bba\label{eqnew1} \sum_{w\in \Ss_n}q^{-l(w)}T_{w^{-1}}\ot \tr^n(T_w) = \sum_{w\in\Ss_{(n-1,1)}}q^{-l(w)}(L_n-[-r]_q)T_{w^{-1}}\ot T_w.\eea\end{lem}
The proof will be given in the Appendix \ref{proofof53}.

 We check \rref{eq23}. By definition of $\etr$,
\bbs {R_w}^{J_1m}_{K_1p}C^p_m=\etr(R_w)^{J_1}_{K_1}.\ees
Thus, without indices, \rref{eq23} has the form
 \bbs 
\sum_{w\in\Ss_n}q^{-l(w)}(\ove{P_n}^{-1}R_{w^{-1}}C^{\ot n-1})\ot \etr^n(R_w)=
\sum_{w\in\Ss_{(n-1,1)}}q^{-l(w)}(\ove{P_{n-1}}^{-1}R_{w^{-1}}C^{\ot n-1})\ot R_w.\ees
Remark that $P_n=P_{n-1}(L_n-[-r])$, and $\ove{P_{n-1}}$ acts only on the first $n-1$ components of $\Vn$, we see that the above equation follows from \rref{eqnew1}.

We check \rref{eq22}. \rref{eq22} can be rewritten as
\bbs \sum_{w\in\Ss_n}q^{-l(w)}((\ove{L_n}-[-r])^{-1}R_{w^{-1}}C^{\ot n-1})_{I_1m}^{L_1m} {R_w}_{K}^J=
\sum_{w\in\Ss_{n-1}n}q^{-l(w)}(R_{w^{-1}}C^{\ot n-1})_{I_1}^{L_1}{R_w}_{K_1}^{J_1}\delta_{k_n}^{j_n}.\ees
Omitting indices and canceling $C^{\ot n-1}$, it has the from 
\bbs 
\sum_{w\in\Ss_n}q^{-l(w)}\etr^n((\Ln{n}-[-r])^{-1}R_{w^{-1}})\ot R_w=\sum_{w\in\Ss_{n-1,1}}q^{-l(w)}R_{w^{-1}}\ot R_w.\ees

The trick here is to take $(\Ln{n}-[-r])^{-1}$ out of $\etr^n$. To do this we consider the left multiplication by this element as an endomorphism of $\ove{\H_n}$.

Let $P\in \H_n$ such that $\ove{P}=(\Ln{n}-[-r]_q)^{-1}$. Let
\bbas(L_n-[-r]_q)T_v=\sum_{u\in\Ss_n}c^u_vT_u \mbox{ and }
PT_v=\sum_{u\in\Ss_n}d^u_vT_u,\eeas
then 
\bbs\sum_{w\in\Ss_n} c^w_v d_w^u R_u=\delta_{u,v}R_v.\ees

Set $<T_u,T_v>=\delta_{u,v}q^{l(v)}$ and extend it linearly on $\H_n$, we get a non-degenerate symmetric scalar product on $\H_n$ which satisfies $<hk,g>=<h,gk^*>$ and $<h,gk>=<g^*h,k>$, where $*$ is the linear extension of the map $T_w\loma T_{w^{-1}}, w\in \Ss_n$ (\cite{dj1}, Lemma 2.2). Since $(L_n-[-r]_q)$ is $*$-invariant, we have:
\bbs q^{l(w)}c_v^w=<T_w,(L_n-[-r]_q)T_v>=<(L_n-[-r]_q)T_w,T_v>.\ees
Hence 
\bba \label{eqnew3}(L_n-[-r]_q)T_w=\sum_{v\in\Ss_n}q^{l(w)-l(v)}c^w_vT_v.\eea
Since $T_w,w\in\Ss_n$ form a basis for $\H_n$, we can conclude that $c^u_v=q^{l(v)-l(u)}c_u^v$. Applying first operation $*$ on \rref{eqnew3} and then $\rho$, we get
\bbs R_{w^{-1}}(L_n-[-r]_q)=\sum_{v\in\Ss_n}q^{l(w)-l(v)}c^w_vR_{v^{-1}}.\ees

Now, we have
\bbas\lefteqn{ \sum_{w\in\Ss_n}q^{-l(w)}\etr^n((\Ln{n}-[-r]_q)^{-1}R_{w^{-1}})\ot R_w(\Ln{n}-[-r]_q)}&&\\ &&=
\sum_{w,u,v\in\Ss_n}q^{-l(w)}d^{v^{-1}}_{w^{-1}}\etr^n(R_{v^{-1}})\ot q^{l(w)-l(u)}c_{u^{-1}}^{w^{-1}}R_u\\ &&=
\sum_{w,u,v\in\Ss_n}q^{-l(u)}d^{v^{-1}}_{w^{-1}}\etr^n(R_{v^{-1}})c_{u^{-1}}^{w^{-1}}\ot R_u\\ &&=
\sum_{u\in\Ss_n}q^{-l(u)}\etr^n(R_{u^{-1}})\ot R_u\mbox{ (by virtue of Lemma \ref{technicallem}) }\\ 
 &&=
\sum_{u\in\Ss_n}q^{-l(v)}(R_{v^{-1}})\ot R_v(\Ln{n}-[-r]_q).\eeas
 Since $(\Ln{n}-[-r]_q)$ is invertible, \rref{eq22} is proved.
\eee
\subsection{The Integral on $SH_R$}
 The integral on $SH_R$ is easier to compute. By its definition, $SH_R$ is a quotion of the bialgebra $E_R$, generated by $\{x_i^j|1\leq i,j\leq d\}$, by setting the quantum determinant equal to 1. Thus, we have to compute $\int(Z_I^J)$, where $I,J$ are multi-indices. 

Introduce a $\bZ_r$-grading on $SH_R$, setting deg$(z_i^j)=1\mod r$. Then $SH_R$ is defined by $\bZ_r$-homogeneous relation. From its definition we see that 
$$\int(Z_I^J)=0,\quad \mbox{unless } l(I)=l(J)=0 \mod r.$$
\begin{thm}\label{thmintshr} Let $l(I)=l(J)=n=kr$ and $F_{(k^r)}$ be the minimal central idempotent, that corresponds to $(k^r)$. Let $\Phi_k=\rho_n(F_{(k^r)})$. Then
\bba\label{intshr} \int(Z_I^J)={\Phi_k}_I^J.\eea\end{thm}
\proof It is known that the integral has the following property. For a comodule $M$ on $SH_R$, the map $M\lora M$, $m\loma m_0\int(m_1)$ is the projection from $M$ onto its subspace of $SH_R$-coinvariants, that are $n\in M, \delta(n)=n\ot 1$. Moreover, the integral is determined by this property.

Since every simple $SH_R$-comodule is a submodule of $\Vn$ for some $n$, it is sufficient to check the above property for $\Vn$. Notice that the subspace of $SH_R$-coinvariant in $\Vn$ is the sum of all $\Vn$-subcomodules that are isomorphic to $\bK$. From the proof of Theorem \ref{thmshr}, each submodule of $\Vn$, isomorphic to $\bK$, is determined by a primitive idempotent of $\H_n$, that corresponds to ${(k^r)}$. In other words, if $\Vn$ contains a non-trivial subspace of $SH_R$-coinvariant, then $n=kr$ for some $k$ and for the image of primitive idempotent, corresponding to $(k^r)$, under $\rho_n$ is non-zero (Corollary \ref{corsimple}). As an easy consequence, we see that the subspace of $SH_R$-coinvariant is precisely $\Im(\rho_n(F_{(k^r)})).$ \eee

In Theorem \ref{dldj}, a full system of mutually orthogonal primitive idempotent, corresponding to a partition $\lambda$, is given. Consequently, their sum is the minimal central idempotent corresponding to $\lambda$:
$$F_\lambda=\sum_{i=1}^{d_\lambda}E_{{i,\lambda}}.$$

We give below a simpler description of $\Phi_k$. First we remark
\begin{lem}\label{lemcentralize} The minimal central idempotent $F_\lambda$ can be obtained from a primitive central idempotent $E_\lambda$ in the following way:
$$ F_\lambda=\tr(E_\lambda)\sum_{w\in\Ss_n}q^{-l(w)}T_wE_\lambda T_{w^{-1}}$$
where $\tr(h):=<h,1>,$ $<,>$ is the scalar product introduced in the previous subsection.\end{lem}
\proof Let $A_\lambda$ be the two-sided ideal generated by $F_\lambda$. Thus, $A_\lambda$ is a full matrix ring, of degree $d_\lambda$. Let $\{E_\lambda^{i,j}|1\leq i,j\leq d_\lambda\}$ be a basis of $A_\lambda$ with the property
$$E_\lambda^{ij}E_\lambda^{kl}=\delta^j_kE_\lambda^{il}.$$ Thus, $E_\lambda^{ii}$ are primitive idempotents and $F_\lambda=\sum_iE_\lambda^{ii}$. We let $\lambda$ run it the set of all partitions of $n$, then $\{E_\lambda^{ij}|\lambda\part n,1\leq i,j\leq d_\lambda\}$ form a basis for $\H_n$. 

The map $\tr$ is by definition a trace map, which is faithful, i.e. $\forall G, \exists H, \tr(GH)\neq 0$.
By standard argument (see \cite{g-h-j}) $\tr(E_\lambda^{ij})=0$ unless $i=j$ and $\tr(E_\lambda^{ii})$ are equal for all $i$ to a non-zero constant $k_\lambda$. Hence $\{\frac{1}{k_\lambda}E_\lambda^{ji}|\lambda\part n,1\leq i,j\leq d_\lambda\}$ form a basis, dual to the described above basis with respect to the scalar product $<,>$. On the other hand, the bases $\{T_w|w\in\Ss_n\}$ and $\{q^{-l(w)}T_w|w\in\Ss_n\}$ are dual w.r.t this scalar product, too. Therefore we have:
\bbs \sum_{w\in\Ss_n}q^{-l(w)}T_wHT_{w^{-1}}=\sum_{\lambda\part n\atop 1\leq i,j\leq d_\lambda}\frac{1}{k_\lambda}E^{ij}_\lambda HE_\lambda^{ji}, \quad \mbox{for any }H\in\H_n.\ees
In particular, for $H=E_\lambda=E_\lambda^{kk}$, we have
\bbs \sum_{w\in\Ss_n}q^{-l(w)}T_wE_\lambda T_{w^{-1}}=\frac{1}{k_\lambda}F_\lambda. \eqeee\ees

The value of $\tr(E_\lambda)$ can be computed explicitly:
\bba\label{eqtrel}\tr(E_\lambda)=q^{\sum_i{\lambda_i}(i-1)}\prod_{k=1}^n\frac{1}{[r_\lambda(k)+r]_q}\prod\frac{[\lambda_i-\lambda_j+j-i]_q}{[j-i]_q}.\eea
The proof of this formula will be given in a forthcoming paper \cite{ph98}.

In our case, the operator $\ove{E_{(k^r)}}$, according to Lemma \ref{lemsimple}, is the product
\bbs \ove{E_{(k^r)}}=\ove{Y_r}\ot \ove{Y_r}\ot\cdots\ot \ove{Y_r}.\ees
Thus, we have
\begin{cor} The operator $\Phi_k$ can be given by
\bba \Phi_k&=&\tr(E_{(k^r)})\sum_{w\in\Ss_{kr}}q^{-l(w)}R_w\ove{Y_r}^{\ot k}R_{w^{-1}}\nonumber\\ 
\label{eqphik}&=&q^{kr(r-1)/2}\frac{[0]_q![1]_q!\cdots[k-1]_q!}{[r]_q![r+1]_q!\cdots[r+k-1]_q!}\sum_{w\in\Ss_{kr}}q^{-l(w)}R_w\ove{Y_r}^{\ot k}R_{w^{-1}}.\eea\end{cor}
\vskip2ex

Comparing the two formulae for integral on $H_R$ and $SH_R$, we notice that the integral on the latter can be obtained from the integral on the other one in the following sense.
	
The quantum determinant can be represented in the form $q^{r(r+1)/2}{B^{\ot r}}_M^P{\ove{Y_r}}^M_NZ^N_P$, and it inverse is $q^{r(r+1)/2}{B^{\ot r}}_M^P{\ove{Y_r}}^M_NT^N_P$, the matrix $B$ was introduced in Equation \rref{eqtz1}. Therefore, if we apply the formula \rref{eqi16} for the element $Z_I^J{B^{\ot r}}_M^P{\ove{Y_r}^{\ot k}}^M_NT^N_P$, $I,J,K,L\in\bN^{kr}$ then we have
\bbas \lefteqn{\int(q^{r(r+1)/2}Z_I^J{B^{\ot kr}}_M^P{\ove{Y_r}^{\ot k}}^M_NT^N_P)}\\
&=&q^{r(r+1)/2} \sum_{w\in\Ss_{kr}}q^{-l(w)}\left(\prod_{m=1}^{kr}(R_{w^{-1}}\Ln{m}-[-r]_q)^{-1}C^{\ot kr}B^{\ot kr} E_{(k^r)}R_w\right)_I^J\\
   &=&q^{kr(r-1)/2}\frac{[0]_q![1]_q!\cdots[k-1]_q!}{[r]_q![r+1]_q!\cdots[r+k-1]_q!}\sum_{w\in\Ss_n}q^{-l(w)}\left(R_w \ove{Y_r}^{\ot k}R_{w^{-1}}\right)_I^J ,\eeas
in the second equation we use Equations \rref{eqap3} and \rref{eqap31} and the fact that $BC=q^{-(r+1) }I$ which follows from Equation \rref{eqinj1}.

\begin{appendix}
\newtheorem{appenv}{}[section]

\section{The Hecke Algebras}
We recall in this Appendix the definition and some important properties of the Hecke algebra. Then we derive from them some results that are needed in our context.

The symmetric groups $\Ss_n$ consists of permutations of the set $\{1,2,\ldots,n\}$. It is generated by the basic transpositions $v_i=(i,i+1), 1\leq i\leq n-1$. The length of a permutation $v$ is the number of pairs $(i,j),1\leq i<j\leq n$, such that $iv>jv$. This is equal to the minimal number of $v_i$'s needed to express $v$. 
Irreducible representations of $\Ss_n$ are indexed by partitions of $n$. 

Let $\lambda$ be a partition of $n$, a $\lambda$-tableau $a^\lambda$ is a matrix $a_i^j,1\leq j\leq \lambda_i,i=1,2,\ldots$, where $a_i^j$ are different elements from the set $\{1,2,\ldots, n\}$. A tableau is said to be row- (column-) standard if the numbers $1,2,\ldots,n$ increase along the rows (columns) and standard if it is row- and column-standard. 

 \begin{appenv}{\bf The Hecke Algebra $\H_n$}\end{appenv}
 Let $\bK$ be a field of characteristic zero, $0\neq q\in k$. As a $\bK$-vector space, the Hecke algebra $\H_{n,q}$ is spanned by $\{T_v|v\in \Ss_n\}$. The multiplication satisfies the following relations:
 \begin{itemize}
 \item[] $T_\unit=\unit_{\H_{n,q}}$,
 \item[] $T_vT_w =T_{vw}$ iff $l(wv)=l(v)+l(w)$,
 \item[] $T_{v_i^2}=(q-1)T_{v_i}+q\unit$, $v_i=(i,i+1)$.
 \end{itemize}
 Let us denote $T_i:=T_{v_i}$, $1\leq i\leq n-1$. As is the case of $\Ss_n$, $T_i$ generate $\H_n$ and satisfy:
 \bbas T_iT_{i+1}T_i=T_{i+1}T_iT_{i+1}, 1\leq i\leq n-2.\eeas


If $q$ is not a root of unity then $\H_n=\H_{n,q}$ is semi-simple.  

For $n\in\bZ$, define $ [n]_q:=(q^n-1)/(q-1)$ and $[\pm n]_q!:=[\pm 1]_q[\pm 2]_q\cdots [\pm n]_q$, $[0]_q!:=1$. We denote,
\bbs  X_n:= \frac{1}{[n]_q!}\sum_{w\in\Ss_{n}} T_w,\quad
Y_n:=\frac{1}{[n]_{1/q}!}\sum_{w\in \Ss_{n}}(-q)^{-l(w)}T_w.\ees
$X_n$ (resp. $Y_n$) is a minimal central primitive idempotent of $\H_n$, which induces the trivial (resp. signature) representation of $\H_n$
\bbs T_wX_n=X_nT_w=q^{l(w)}X_n,\quad T_wY_n=Y_nT_w=(-1)^{l(w)}Y_n.\ees
\begin{appenv}{Murphy Operators and Primitive Idempotents}\end{appenv}
Primitive idempotents of $\H_n$ can be expressed in terms of Murphy operators. The Murphy operators $L_m$, $1\leq m\leq n$, were introduced in \cite{dj2} as follows: $L_1:=0,$
\bbas 
 L_m:= q^{-1}T_{m-1}+q^{-2}T_{m-2}T_{m-1}T_{m-2}+\cdots+q^{-m+1}T_1T_2\cdots T_{m-1}\cdots T_1,m\geq 2.\eeas
The following equality will be frequently used:
\bbas T_kT_{k+1}\cdots T_{m-1}T_mT_{m-1}\cdots T_{k+1}T_k=T_mT_{m-1}\cdots T_{k+1}T_kT_{k+1}\cdots T_{m-1}T_m, \\ (k< m\leq n-1).\eeas
We summarize some properties of Murphy operators obtained in \cite{dj2} in a theorem.
\begin{appenv}{\bf Theorem}\label{dldj} \it
\begin{itemize}
\item[(i)] $L_m$, $1\leq m\leq n$ commute and symmetric polynomials on them form the center of $\H_n$ (\cite{dj2}, Theorem 2.14).
			       	      
\item[(ii)] Simple $\H_n$ modules are parameterized by partitions of $n$. The simple module $S_\lam$, $\lam\part n$, has dimension $d_\lambda$ equal to the number of standard $\lambda$-tableaux. There exists a basis (called Young's semi-normal basis), indexed by standard $\lam$-tableaux, $l_i=l_{a_i^\lambda}, i=1,2,\ldots,d_\lambda$, such that 
\begin{equation}\label{eqap2}l_iL_m=[r_{i,\lam}(m)]_ql_i,\end{equation}		          
$r_{i,\lam}(m):=k-j,$ where $(j,k)$ is the coordinate of $m$ in the standard $\lam$-tableau $a_i^\lambda$ (\cite{dj2}, Lemma 4.6). \footnote{ In \cite{dj2}, $r_{i,\lam}(m)$ is defined to be $[k-j]_q$.}
\item[(iii)] The set
\bbs \left\{ \left. E_{i,\lam}:=\prod_{|k|\leq n-1}\prod_{1\leq m\leq n\atop k\neq r_{i,\lam}(m)} \frac{L_m-[k]_q}{[r_{i,\lam}(m)]_q-[k]_q}\right| 1\leq i\leq d, \lam\part n \right\}\ees      
is a complete set of primitive idempotents of $\H_n$ (\cite{dj2}, Theorem 5.2).\end{itemize}\end{appenv}

As a direct consequence of Theorem \ref{dldj} we have
\begin{equation}\label{eqap3}\prod_{|c|\leq m-1}(L_m-[c]_q)=0,\end{equation}
since $|r_{{i,\lambda}}(m)|\leq m-1$. Hence we have
\bba \label{eqap31}
E_\lambda=\prod_{{1\leq m\leq n\atop |k|\leq m-1}\atop k\neq r_{i,\lam}(m)}\frac{{L_m}-[k]_q}{[r_{i,\lam}(m)]_q-[k]_q}.\eea

\begin{appenv}{\bf Corollary}\it\label{cortwn} For any primitive idempotent $E\in\H_n$, we have
\bbs T_{w_n}^2E_\lambda=q^{-s(\lambda)-n(n-1)/2}E_\lambda.\ees
\end{appenv}
\proof According to \rref{eqtau2},
$T_{w_n}^2E=\prod_{k=1}^nq^k(1+(q-1)L_k).$
Equation \rref{eqap3} and Formula \rref{eqap31} imply the desired equation.\eee

\begin{appenv}{\bf The completion of  Lemma \ref{lemsimple}}\label{appprim}\end{appenv}
 Assume that $n=m+r$, $r> 1$. For an element $H\in\H_n$ denote $\overline{H}:=\rho(h)$ and for an element $K\in\H_m$, $K^r$ its image under the embedding $\H_m\ni T_j\lora T_{r+j}\in \H_n$. Assume that
\bbs E=E_\lambda=\prod_{{1\leq k\leq m\atop |c|\leq k-1}\atop c\neq r_{i,\lambda}(k)}\frac{\ove{L^r_k}-[c]_q}{[r_{i,\lam}(k)]_q-[c]_q},\ees
where $\lam$ is a partition of $m$.

We have
$\ove{Y_r}(\ove{1}-\ove{L_{r+1}})=[r+1]_q\ove{Y_{r+1}}=0.$
On the other hand, for $k\geq 1$, 
\bbs L_{k+r}=L_k^r+ q^{-k+1}T_{k+r-1}\cdots T_{r+1}L_{r+1}T_{r+1}\cdots T_{k+r-1}.\ees
Hence
\bbs \ove{Y_{r}}\ove{L_{k+r}}= \ove{y_{r}}\ove{L_k^r}+q^{-k+1} \ove{T_{k+r-1}}\cdots \ove{T_{r+1}}\ove{Y_r}\ove{L_{r+1}}\ove{T_{r+1}}\cdots \ove{T_{k+r-1}} 
=\ove{Y_r}(1+q\ove{L_k^r}).\ees
Therefore
\begin{equation}\label{eqstar}
\ove{Y_r}(\ove{L_{k+r}}-[c]_q)=\ove{Y_r}(1+q\ove{L_k^r}-[c]_q)
=q\ove{Y_r}(\ove{L_k^r}-[c-1]_q).\end{equation}

Remark that $\ove{Y_r}$ commutes with all $\ove{L_k^r}$, thus 
\bbs \overline{Y_rE^r}=\ove{Y_r}\prod_{1\leq k\leq m\atop {|c|\leq k-1\atop c\neq r_{i,\lam}(k)}}\frac{q^{-1}(\ove{L_{k+r}}-[c+1]_q)}{[r_{i,\lam(k)}]_q-[c]_q}= \ove{Y_r}\prod_{1\leq k\leq m\atop {|c|\leq k-1\atop c\neq r_{i,\lam}(k)}}\frac{\ove{L_{k+r}}-[c+1]_q}{[r_{i,\lam(k)}+1]_q-[c+1]_q}  .\ees
Equations \rref{eqstar} and \rref{eqap3} imply
\bbs\ove{Y_r}\prod_{|c|\leq k-1}(\ove{L_{k+r}}-[c+1]_q)=\ove{Y_r}\prod_{|c|\leq k-1}(\ove{L^r_k}-[c+1]_q)=0,\ees
which means that
\bbs\ove{Y_r}\prod_{|c|\leq k-1\atop c\neq r_{i,\lam}(k)}\frac{\ove{L_{k+r}}-[c+1]_q}{[r_{i,\lam(k)}+1]_q-[c+1]_q}\ \cdot \ \ove{L_{k+r}}=
[r_{i,\lam(k)}+1]\ove{Y_r}\prod_{|c|\leq k-1\atop c\neq r_{i,\lam}(k)}\frac{\ove{L_{k+r}}-[c+1]_q}{[r_{i,\lam(k)}+1]_q-[c+1]_q}.\ees
Hence, for $r\geq 2$, $k=1,2,\ldots,m$,
\bbs\ove{Y_r}\prod_{|c+1|\leq k+r-1\atop c\neq r_{i,\lam}(k)}\frac{\ove{L_{k+r}}-[c+1]_q}{[r_{i,\lam(k)}+1]_q-[c+1]_q}=
\ove{Y_r}\prod_{|c|\leq k-1\atop c\neq r_{i,\lam}(k)}\frac{\ove{L_{k+r}}-[c+1]_q}{[r_{i,\lam(k)}+1]_q-[c+1]_q}.\ees
Since $\ove{L_{k+r}}$ commute each other, $Y_r$ is an idempotent, we have:
\bbs 
\ove{Y_r}\prod_{1\leq k\leq m\atop {|c+1|\leq k+r-1\atop c\neq r_{i,\lam}(k)}}\frac{\ove{L_{k+r}}-[c+1]_q}{[r_{i,\lam(k)}+1]_q-[c+1]_q}=\ove{Y_rE^r}.\ees
The left-hand side term of this equality is a primitive idempotent, hence, so is the right-hand side term.\eee

\begin{appenv}{\bf Lemma} \label{lemtrq}  Let $R$ be an even Hecke symmetry of rank $r$ and $\rho_n:H_n\lora \ove{\H_n}$ be the corresponding representation of $\H_n$. Then $(\Ln{n}-[-r]_q)$ is invertible in $\ove{\H_n}$.\end{appenv}
\proof  According to Proposition \ref{pror+1}, $\ove{Y_\lambda}=0$ whenever $\lambda_{r+1}\neq 0$. Hence a simple module $S_\lambda$ of $\H_n$ remains simple $\ove{\H_n}$-module only if $\lambda_{r+1}=0$. For these $\lambda$, $r_{{i,\lambda}}(m)\geq -r+1$. Hence, according to Theorem \ref{dldj}, (ii), $\Ln{n}-[-r]_q$ is invertible on $S_\lambda$. This holds for all $\lambda$, thus, $\Ln{n}-[-r]_q$ is invertible on $\ove{\H_n}$.\eee

\begin{appenv}{\bf The proof of Lemma \ref{technicallem}}\label{proofof53}\end{appenv}
the left hand side of Equation \rref{eqnew1} is equal to
\bbs -[-r]\sum_{w\in\Ss_{(n-1,1)}}q^{-l(w)}T_{w^{-1}}\ot T_w+
\sum_{k=1}^{n-1}\sum_{w\in\Ss_{(n-1,1)}}q^{-l(w)-k}T_{n-k}\cdots T_{n-1}T_{w^{-1}}\ot T_wT_{n-2}\cdots T_{n-k}.\ees

We have
\bbas
\lefteqn{\sum_{w\in\Ss_{(n-1,1)}}q^{-l(w)-k}T_{n-k}\cdots T_{n-1}T_{w^{-1}}\ot T_wT_{n-2}\cdots T_{n-k}}&&\\
&& =
\sum_{w\in \Ss_{(n-1,1)}\atop l(wv_{n-2})=l(w)+1}\left[
\bbar{r}q^{-l(w)-k}T_{n-k}\cdots T_{n-2} T_{n-1}T_{w^{-1}}\ot T_{wv_{n-2}}T_{n-3}\cdots T_{n-k}\\
 +q^{-l(w)-k-1}T_{n-k}\cdots T_{n-2}T_{n-1} T_{n-2}T_{w^{-1}}
\ot T_w{(T_{n-2})}^2 T_{n-3}\cdots T_{n-k}\eear\right]\hspace{30ex}\eeas\vskip-1ex
\bbas &&=
\sum_{w\in \Ss_{(n-1,1)}\atop l(wv_{n-2})=l(w)+1}q^{-l(w)-k-1}\left[\bbar{r}
qT_{n-k}\cdots T_{n-2}T_{n-1} T_{w^{-1}}\ot T_{wv_{n-2}}T_{n-3}\cdots T_{n-k}\\
+qT_{n-1}T_{n-k}\cdots T_{n-2}T_{n-1}T_{w^{-1}} \ot T_{w} T_{n-3}\cdots T_{n-2}\\
+(q-1)T_{n-k}\cdots T_{n-2}T_{w^{-1}} \ot T_{wv_{n-2}}  T_{n-3}\cdots T_{n-k}
\eear\right]\hspace{30ex}\eeas\vskip-1ex
\bbas  &&=
\sum_{w\in \Ss_{(n-1,1)}\atop l(wv_{n-2})=l(w)+1}q^{-l(w)-k-1}\left[\bbar{r}
T_{n-k}\cdots T_{n-1}{(T_{n-2})}^2T_{w^{-1}} \ot T_{wv_{n-2}} T_{n-3}\cdots T_{n-2}\\
+qT_{n-1}T_{n-k}\cdots T_{n-2}T_{n-1}T_{w^{-1}}\ot T_wT_{n-3}\cdots T_{n-k}
\eear\right]\hspace{30ex}\eeas\vskip-1ex
\bbas  &&=
\sum_{w\in \Ss_{(n-1,1)}\atop l(wv_{n-2})=l(w)+1}\left[\bbar{r}
q^{-l(wv_{n-2})-k}T_{n-1}T_{n-k}\cdots T_{n-2}T_{n-1}T_{wv_{n-2}}^{-1} \ot T_{wv_{n-2}}T_{n-3}\cdots T_{n-k}\\
+q^{-l(w)-k}T_{n-1}T_{n-k}\cdots T_{n-1}T_{w^{-1}}\ot T_wT_{n-3}\cdots T_{n-k}
\eear\right]\hspace{30ex}\eeas\vskip-1ex
\bbas &&=
\sum_{w\in\Ss_{(n-1,1)}}q^{-l(w)-k}T_{n-1}T_{n-k}\cdots T_{n-2}T_{n-1}T_{w^{-1}}\ot T_wT_{n-3}\cdots T_{n-k}\hspace{30ex}\eeas\vskip-1ex
  \bbas &&=
\cdots \mbox{ (repeating the above process $k-2$ times)\hspace{40ex} }\hspace{30ex}\hspace{30ex} \eeas\vskip-1ex
\bbas &&= 
\sum_{w\in\Ss_{(n-1,1)}}q^{-l(w)-k}T_{n-1}T_{n-2}\cdots T_{n-k}T_{n-k+1}\ldots T_{n-2}T_{n-1}T_{w^{-1}}\ot T_w
\hspace{50ex}\eeas
Summing up the above equations for $k=1,2,\ldots n-1$ we arrive at the desired equation.

\end{appendix}


\begin{thebibliography}{10}

\bibitem{dj1}
R.~{D}ipper and G.~James.
\newblock {Representations of {H}ecke Algebras of {G}eneral Linear Groups}.
\newblock {\em Proc. London Math. Soc.}, 52(3):20--52, 1986.

\bibitem{dj2}
R.~{D}ipper and G.~James.
\newblock {Block and Idempotents of {H}ecke Algebras of {G}eneral Linear
  Groups}.
\newblock {\em Proc. London Math. Soc.}, 54(3):57--82, 1987.

\bibitem{dpw}
J.~{D}u, B.~Parshall, and J.~Wang.
\newblock {Two-parameter Quantum Linear Groups and the Hyperbolic Invariance or
  $q$-{S}chur Algebras}.
\newblock {\em Journal of London Math. Soc.}, 44(2):420--436, 1991.

\bibitem{c-g}
J.-L.Gervais E.Cremer.
\newblock {The Quantum Groups Structure Associated With Non-linearly Extended
  Virasoro Algebras}.
\newblock {\em Comm. Math. Phys.}, 134:619--632, 1990.

\bibitem{frt}
L.~{F}addeev, N.~Reshetikhin, and L.~Takhtajian.
\newblock {Quantisation of {L}ie Groups and {L}ie Algebras}.
\newblock {\em Leningrad Math. Journal}, 1:193--225, 1990.

\bibitem{g-h-j}
F.M. Goodman, P.~de~la Harpe, and V.F.R. Jones.
\newblock {\em {Coxter Graphs and Towers of Algebras}}, volume~14 of {\em MSRI
  publications}.
\newblock Springer-Verlag, 1989.

\bibitem{green2}
J.A. Green.
\newblock {Locally Finite Representations}.
\newblock {\em Journal of Algebra}, 41:137--171, 1976.

\bibitem{gur}
D.I. {G}urevich.
\newblock {Hecke Symmetry and Quantum Determinants}.
\newblock {\em Soviet Math. Doklad}, 38:555--559, 1989.
\newblock See also {\em Leningrad Math. Journal}, 2(4):801-828,1991.

\bibitem{gur1}
D.I. Gurevich.
\newblock { Algebraic Aspects of the Quantum Yang-Baxter Equation}.
\newblock {\em Leningrad Math. Journal}, 2(4):801--828, 1991.

\bibitem{gur96}
D.I. Gurevich, P.N. Pyatov, and P.A. Saponov.
\newblock {Hecke Symmetries and Characteristic Relations on Reflection Equation
  Algebras}.
\newblock {\em q-alg}, 9605048, 1996.

\bibitem{ph97}
Phung~Ho Hai.
\newblock {{K}oszul Property and {P}oincar\'{e} Series of Matrix Bialgebra of
  Type ${A}_n$}.
\newblock {\em Journal of Algebra}, 192(2):734--748, 1997.

\bibitem{ph98}
Phung~Ho Hai.
\newblock {Characters for Quantum Groups of Type $A_n$}.
\newblock {\em Preprint ICTP, {\tt xxx.lanl.gov/dvi/q-alg/98-7045}}, 1998.

\bibitem{haya1}
T.~{H}ayashi.
\newblock { Quantum Groups and Quantum Determinant}.
\newblock {\em Publ. RIMS, Kyoto Univ.}, 28:57--81, 1992.

\bibitem{james1}
G.D. James.
\newblock {\em {The Representation Theory of Symmetric Groups}}.
\newblock Lecture Note in Mathematics. Springer Verlang, Berlin, 1978.

\bibitem{js3}
A.~Joyal and R.~Street.
\newblock {Tortitle Yang-Baxter Operators in Tensor Categories}.
\newblock {\em J. Pure Appl. Algebra}, pages 43--51, 1991.
\newblock MR 92e:18006.

\bibitem{kt1}
C.~Kassel and V.G. Turaev.
\newblock {Double Construction for Monoidal Categories}.
\newblock {\em Acta Mathematica}, 175(1):1--48, 1995.

\bibitem{kassel}
Ch~Kassel.
\newblock {\em Quantum Groups}, volume 155 of {\em Graduate Texts in
  Mathematics}.
\newblock Springer-Verlag, 1995.
\newblock 531p.

\bibitem{k-w}
D.~Kazhdan and H.~Wenzl.
\newblock Reconstructing monoidal categories.
\newblock In Sergej (ed.) et~al. Gelfand, editor, {\em I. M. Gelfand seminar.
  Part 2: Papers of the Gelfand seminar in functional analysis held at Moscow
  University}, volume~16 of {\em Adv. Sov. Math.}, pages 111--136, Providence,
  RI: American Mathematical Society, 1993.

\bibitem{l-t}
R.~{L}arson and J.~Towber.
\newblock { Two Dual Classes of Bialgebras Related To The Concepts of
  ``{Q}uantum Groups'' and ``{Q}uantum {L}ie Algebra''}.
\newblock {\em Comm. in Algebra}, 19(12):3295--3345, 1991.

\bibitem{lyu1}
V.V. Lyubashenko.
\newblock { Hopf Algebras and Vector Symmetries}.
\newblock {\em Russian Math. Survey}, 41(5):153--154, 1986.

\bibitem{ls}
V.V. {L}yubashenko and A.~Sudbery.
\newblock {Quantum Super Groups of {GL}$(n|m)$ Type: Differential Forms,
  {K}oszul Complexes and {B}erezinians}.
\newblock {\em Duke Math. Journal}, 90:1--62, 1997.

\bibitem{mcd2}
I.G. Macdonald.
\newblock {\em {Symmetric functions and the Hall polynomials}}.
\newblock Oxford University Press, 1979 (Second edition 1995).
\newblock New York.

\bibitem{majid3}
S.~{M}ajid.
\newblock {Quasitriangular {H}opf Algebras and {Y}ang-{B}axter Equation}.
\newblock {\em International Journal of Modern Physics A}, 5:1--91, 1990.

\bibitem{manin1}
Yu.I. {M}anin.
\newblock {\em {Quantum Groups and {N}on-commutative {G}eometry}}.
\newblock GRM, Univ. de Montreal, 1988.

\bibitem{manin2}
Yu.I. {M}anin.
\newblock {Multiparametric Quantum Deformation of the {G}eneral Linear
  Supergroups}.
\newblock {\em Com. Math. Phys.}, 123:163--175, 1989.

\bibitem{ps96}
P.N. Pyatov and P.A. Saponov.
\newblock {Newton Relations for Quantum Matrix Algebras of ${RTT}$-Type}.
\newblock {\em Preprint IHEP,96-76}, 1996.

\bibitem{rt}
N.~Reshetikhin and V.~Turaev.
\newblock {Ribbon Graph and Their Invariant Derived from Quantum Groups}.
\newblock {\em Comm. Math. Phys.}, 127:1--26, 1990.

\bibitem{schauen}
P.~{S}chauenburg.
\newblock {On Coquasitriangular {H}opf Algebras and the Quantum {Y}ang-{B}axter
  Equation}.
\newblock {\em Algebra Berichte 67, Verlag Reinhard Fischer, Munich}, 1992.

\bibitem{takeuchi}
M.~Takeuchi.
\newblock {Free Hopf Algebra Generated by Coalgebras}.
\newblock {\em Journal Math. Soc. Japan}, 23:561--582, 1971.

\bibitem{tur1}
V.G. Turaev.
\newblock {Modular Categories and 3-Manifolds Invariant}.
\newblock {\em Intern. J. Modern Phys. B}, 6:1807--1824, 1992.

\end{thebibliography}
\end{document}